\documentclass[11pt,a4paper]{article}

\usepackage{amssymb}
\usepackage{amsmath}
\usepackage{epic}
\usepackage{amsthm}
\usepackage{amscd}

\newtheorem{Thm}{Theorem}[section]
\newtheorem{Cor}[Thm]{Corollary}
\newtheorem{Lem}[Thm]{Lemma}
\newtheorem{Prop}[Thm]{Proposition}
\newtheorem{Def}[Thm]{Definition}

\theoremstyle{remark}
\newtheorem*{remark}{\rm\textbf{Remark}}

\newtheorem*{pf}{\rm\textbf{Proof}}

\newtheorem*{ac}{\rm\textbf{Acknowledgement}}

\begin{document}

\title{\Large \bf On the Affine Schur Algebra of Type $A$ }

\author{\large Dong Yang\footnote{The author acknowledges support by the AsiaLink network
Algebras and Representations in China and Europe,
ASI/B7-301/98/679-11.}}

\date{}

\maketitle \begin{quote} \abstract{ Let $n,r\in\mathbb{N}$. The
affine Schur algebra $\widetilde{S}(n,r)$ (of type A) over a field
$K$ is defined to be the endomorphism algebra of certain tensor
space over the extended affine Weyl group of type $A_{r-1}$. By the
affine Schur--Weyl duality it is isomorphic to the image of the
representation map of the $\mathcal{U}(\widehat{\mathfrak{gl}}_{n})$
action on the tensor space when $K$ is the field of complex numbers.
We show that $\widetilde{S}(n,r)$ can be defined in another two
equivalent ways. Namely, it is the image of the representation map
of the semigroup algebra $K\widetilde{GL}_{n,a}$ (defined in
Section~\ref{S:semigroups}) action on the tensor space and it equals
to the 'dual' of a certain formal coalgebra related to this
semigroup. By these approaches we can show many relations between
different Schur algebras
 and affine Schur algebras and reprove one side of the affine Schur--Weyl duality.}
\end{quote}

\section{Introduction}\label{S:introduction}

The notion of the affine $q$-Schur algebra was first introduced by
R.M.Green~\cite{RMG}, although the algebra was first studied by
V.Ginzburg and E.Vasserot~\cite{GV}. For $n,r\in\mathbb{N}$, the
affine $q$-Schur algebra $\widehat{S}_{q}(n,r)$ is defined as the
endomorphism algebra of the $q$-tensor space (depending on $n$) over
the extended affine Hecke algebra of type $A_{r-1}$. Specializing
this approach at $q=1$ and tensoring the resulting
$\mathbb{Z}$-algebra with an infinite field $K$ we obtain the affine
Schur algebra over $K$ (see Section~\ref{S:affineschuralgebra} for
detailed definition), denoted by $\widetilde{S}(n,r)$.

The purpose of this paper is to relate the affine Schur algebra to
the representation theory of semigroups analogously to the finite
case, and to study the relation between Schur algebras and affine
Schur algebras with the same or different parameters as well as
the relation between 
the affine complex general linear Lie algebra
$\widehat{\mathfrak{gl}}_{n}$ and the affine complex special linear
Lie algebra $\widehat{\mathfrak{sl}}_{n}$ and the affine Schur
algebra over $\mathbb{C}$. We introduce a semigroup
$\widetilde{GL}_{n,a}$ for each $a\in K^{\times}$ and show that it
acts on an infinite dimensional tensor space $\widetilde{E}^{\otimes
r}$. Denote by $\phi$ the corresponding representation map. We
define a formal coalgebra $\widetilde{A}(n,r)$, which is a set of
functions on $\widetilde{GL}_{n,a}$ (for definition of formal
coalgebras see Appendix 1). The main results are the following.

\begin{Thm}\label{T:m1} {\rm (i)}
The affine Schur algebra $\widetilde{S}(n,r)$ is isomorphic to the
image $\rm{Im}$$(\phi)$ of the representation map $\phi$.

{\rm (ii)} The algebra $\widetilde{S}(n,r)$ is isomorphic to
$\widetilde{A}(n,r)^{\#}$ (for a $K$-vector space $V$ with a fixed
basis $\{v_{i}|i\in I\}$, $V^{\#}$ is the sub $K$-vector space of
$V^{*}$ with basis dual to $\{v_{i}|i\in I\}$. See Appendix 1).

\end{Thm}

Theorem~\ref{T:m1}(i)(ii) are analogous to J.A.Green's results for
finite case in~\cite{JAG} \S2.

\begin{Thm}\label{T:m2}
{\rm (i)} There is a natural algebra embedding from the Schur
algebra $S(n,r)$ to the affine Schur algebra $\widetilde{S}(n,r)$.

{\rm (ii)} For each $a\in K^{\times}$ there is a surjective
algebra homomorphism $\psi_{a}:\widetilde{S}(n,r)\rightarrow
S(n,r)$.

{\rm (iii)}  For each $a\in K^{\times}$ there is a surjective
algebra homomorphism
$\widetilde{\rm{det}}_{a}^{\#}:\widetilde{S}(n,n+r)\rightarrow
\widetilde{S}(n,r)$.

{\rm (iv)} These homomorphisms are compatible (for details see
Theorem~\ref{T:gpandalg}).
\end{Thm}

Theorem~\ref{T:m2} (i)(ii)(iv) states the relation between the Schur
algebra and the affine Schur algebra. The quantized version of
Theorem~\ref{T:m2} (i) is \cite{RMG} Proposition 2.2.5. The
quantized version of Theorem~\ref{T:m2} (iii) is studied
by~\cite{Du}~\cite{RMG3} in finite case and by~\cite{L}~\cite{L2} in
affine case.

By~\cite{CP} the affine Schur--Weyl duality holds. Precisely, the
duality says that the quotients of
$\mathcal{U}(\widehat{\mathfrak{gl}}_{n})$ and the group algebra of
the extended affine Weyl group of type $A_{r-1}$ acting faithfully
on the tensor space $\widetilde{E}^{\otimes r}$ centralize each
other. For the quantized version we refer
to~\cite{CP}~\cite{GRV}~\cite{RMG}~\cite{VV}. We can prove one side
of the affine Schur--Weyl duality using the affine Schur algebra.
Precisely, we have the following theorem.

\begin{Thm}\label{T:m3} Assume $K=\mathbb{C}$, $n\geq 2$.

{\rm (i)} There is a surjective algebra homomorphism
$\widetilde{\pi}:\mathcal{U}(\widehat{\mathfrak{gl}}_{n})\rightarrow
\widetilde{S}(n,r)$.

{\rm (ii)} When $r<n$,
$\widetilde{\pi}|_{\mathcal{U}(\widehat{\mathfrak{sl}}_{n})}:\mathcal{U}(\widehat{\mathfrak{sl}}_{n})\rightarrow
\widetilde{S}(n,r)$ is surjective.
\end{Thm}

Recently, the kernel of the quantized $\widetilde{\pi}$ is described
by S.R.Doty and R.M.Green~\cite{DR} in case $r<n$, and by
K.McGerty~\cite{McGerty} in general.

\vskip10pt

This paper is organized as follows. In Section~\ref{S:schuralgebra}
we briefly recall the Schur algebra defined by J.A.Green
in~\cite{JAG} and further developed in~\cite{JAG3} and~\cite{JAG5}
and a great many papers by others. We will generalize most of this
section to the affine Schur algebra in later sections. In
Section~\ref{S:semigroups} we will introduce two semigroups
$\widetilde{GL}_{n,a}$ and $\widetilde{SL}_{n,a}$, which are
analogous to the general and special linear groups. In
Section~\ref{S:affineschuralgebra} we will define the affine Schur
algebra in different but equivalent ways and give two multiplication
formulas in terms of certain basis elements. In
Sections~\ref{S:relations}--\ref{S:liealgebras}, we will give
relations between different Schur algebras and affine Schur algebras
and also between the loop algebras (consequently
$\widehat{\mathfrak{sl}}_{n}$ and $\widehat{\mathfrak{gl}}_{n}$) and
the affine Schur algebra. As a byproduct, we obtain sets of algebra
generators of the affine Schur algebra. We will define and study the
formal coalgebra in Appendix 1. Appendix 2 is a generalization of
J.A.Green's result (Mackey's formula) in~\cite{JAG2} which provides
necessary lemmas for the proof of the second multiplication formula
in Section~\ref{S:affineschuralgebra}.

In the sequel $K$ will be an infinite field. Without further
comment tensor products will be taken over $K$.

\section{The Schur algebra $S(n,r)$}\label{S:schuralgebra}

This section briefly recalls on the equivalent definitions  and
some basic properties of the Schur algebra. The main references
are~\cite{CL}~\cite{JAG}~\cite{JAG3}.

Fix $n\in\mathbb{N}$. Let $\mathcal{M}_{n}=\mathcal{M}_{n}(K)$ be
the algebra of $n\times n$ matrices with entries in $K$,
$GL_{n}=GL_{n}(K)$ the general linear group, and $SL_{n}=SL_{n}(K)$
the special linear group.

Fix $r\in\mathbb{N}$. Let
$I(n,r)=\{\underline{i}=(i_{1},\ldots,i_{r})|i_{t}\in
\{1,\ldots,n\}, t=1,\ldots,r\}$.  The symmetric group $\Sigma_{r}$
on $r$ letters acts on $I(n,r)$ on the right by place permutation,
i.e.
$(i_{1},\ldots,i_{r})\sigma=(i_{\sigma(1)},\ldots,i_{\sigma(r)})$
and on $I(n,r)\times I(n,r)$ diagonally. The equivalence relations
on $I(n,r)$ and on $I(n,r)\times I(n,r)$ induced by this action are
both denoted by the symbol $\sim_{\Sigma_{r}}$.

\subsection{The coordinate ring}

Since the field $K$ is infinite, the coordinate functions $c_{ij}$
on $GL_{n}$, $i,j=1,\ldots,n$ are algebraically independent. So
$A(n)=K[c_{ij}]_{i,j=1,\ldots,n}$ is a polynomial ring in $n^{2}$
indeterminates. Let $c_{\underline{i},\underline{j}}$ denote the
product $c_{i_{1}j_{1}}\cdots c_{i_{r}j_{r}}$. Then
$c_{\underline{i},\underline{j}}=c_{\underline{k},\underline{l}}$ if
and only if
$(\underline{i},\underline{j})\sim_{\Sigma_{r}}(\underline{k},\underline{l})$.
The homogeneous component $A(n,r)$ of $A(n)$ of degree $r$ has a
basis $\{c_{\underline{i},\underline{j}}\}$, where
$(\underline{i},\underline{j})$ runs over a set of representatives
of $\Sigma_{r}$-orbits of $I(n,r)\times I(n,r)$. Moreover, $A(n,r)$
is a coalgebra with respect to the comultiplication $\Delta:
c_{\underline{i},\underline{j}}\mapsto \sum_{\underline{s}\in
I(n,r)}{c_{\underline{i},\underline{s}}\otimes
c_{\underline{s},\underline{j}}}$ and the counit
$\epsilon:A(n,r)\rightarrow K$,
$c_{\underline{i},\underline{j}}\mapsto
\delta_{\underline{i},\underline{j}}$. This makes $A(n)$ into a
bialgebra.

The {\em Schur algebra}, denoted by $S(n,r)$, is defined as the dual
of the coalgebra $A(n,r)$. Let
$\{\xi_{\underline{i},\underline{j}}\}$ be the basis dual to
$\{c_{\underline{i},\underline{j}}\}$, i.e. for
$\underline{i},\underline{j},\underline{p},\underline{q}\in I(n,r)$
\[
\xi_{\underline{i},\underline{j}}(c_{\underline{p},\underline{q}})=\begin{cases}
1, & {\rm
if}\ \ (\underline{i},\underline{j})\sim_{\Sigma_{r}}(\underline{p},\underline{q})\\
0, & {\rm otherwise}.\end{cases}
\]
Then
$\xi_{\underline{i},\underline{j}}=\xi_{\underline{k},\underline{l}}$
if and only if
$(\underline{i},\underline{j})\sim_{\Sigma_{r}}(\underline{k},\underline{l})$.
Multiplication is given by Schur's product rule (\cite{JAG}(2.3b))
\[
\xi_{\underline{i},\underline{j}}\xi_{\underline{k},\underline{l}}=\sum_{(\underline{p},\underline{q})\in
(I(n,r)\times
I(n,r))/\Sigma_{r}}Z(\underline{i},\underline{j},\underline{k},\underline{l},\underline{p},\underline{l})\xi_{\underline{p},\underline{q}}
\]
where
$Z(\underline{i},\underline{j},\underline{k},\underline{l},\underline{p},\underline{l})=\#\{\underline{s}\in
I(n,r)|(\underline{i},\underline{j})\sim_{\Sigma_{r}}(\underline{p},\underline{s}),\
{\rm and\ }
(\underline{s},\underline{q})\sim_{\Sigma_{r}}(\underline{k},\underline{l})
\}$. This product rule is generalized to Coxeter system of type $B$
by R.M.Green~\cite{RMG4}.

One observes that
$\xi_{\underline{i},\underline{j}}\xi_{\underline{k},\underline{l}}=0$
unless $\underline{j}\sim_{\Sigma_{r}}\underline{k}$. J.A.Green
gives another version of the multiplication formula
(\cite{JAG3}(2.6))
\begin{equation}\label{E:greenprod}
\xi_{\underline{i},\underline{j}}\xi_{\underline{j},\underline{l}}=\sum_{\delta\in
\Sigma_{\underline{l},\underline{j}}\backslash\Sigma_{\underline{j}}/\Sigma_{\underline{i},\underline{j}}}[\Sigma_{\underline{i},\underline{l\delta}}:\Sigma_{\underline{i},\underline{j},\underline{l\delta}}]\xi_{\underline{i},\underline{l\delta}}
\end{equation}
where $\Sigma_{\underline{i}}$ is the stabilizer of $\underline{i}$
in $\Sigma_{r}$,
$\Sigma_{\underline{i},\underline{l\delta}}=\Sigma_{\underline{i}}\cap
\Sigma_{\underline{l}\delta}$, and
$\Sigma_{\underline{i},\underline{j},\underline{l\delta}}=\Sigma_{\underline{i}}\cap\Sigma_{\underline{j}}\cap
\Sigma_{\underline{l}\delta}$.

Let $\text{det}$$=\sum_{\sigma\in\Sigma_{n}}{{{\rm
sgn}(\sigma)}c_{(1,\ldots,n),(1,\ldots,n)\sigma}}\in A(n,n)$ where
${\rm sgn}()$ is the sign function. Then multiplication by
$\text{det}$ gives an injective coalgebra homomorphism from $A(n,r)$
to $A(n,n+r)$. Taking the dual we obtain a surjective algebra
homomorphism $\text{det}^{*}: S(n,n+r)\rightarrow S(n,r)$.

\subsection{The representation $E^{\otimes r}$ of
$GL_{n}$}\label{S:naturalrep}

Let $E$ be an $n$-dimensional $K$-vector space. Then $GL_{n}$ acts
on $E$ from the left by matrix multiplication with elements in $E$
considered as column vectors, and diagonally on the $r$-fold tensor
product $E^{\otimes r}$ of $E$. The image of the corresponding
representation map is isomorphic to the Schur algebra $S(n,r)$.

Therefore, there is a surjective map $e: KGL_{n}\rightarrow S(n,r)$.
Actually, $e(g)(c)=c(g)$ for any $c\in A(n,r)$ and $g\in GL_{n}$.
Moreover, the restriction map $e|_{KSL_{n}}$ is surjective and
compatible with $\text{det}^{*}$, i.e. $\text{det}^{*}\circ
e^{n+r}=e^{r}$ on $KSL_{n}$ (see~\cite{Liu}).

There is an equivalence between $S(n,r)$-$\text{mod}$ and $M(n,r)$,
the category of homogeneous polynomial representations of $GL_{n}$
of degree $r$. The space $A(n,r)$ is a homogeneous representation as
well as a homogeneous anti-representation of $GL_{n}$ of degree $r$,
and therefore an $S(n,r)$-bimodule with actions (\cite{JAG} (2.8b))
\[\xi \circ c=\sum\xi(c_{t}') c_{t},~~ c\circ\xi=\sum\xi(c_{t})c_{t}'\]
for $\xi\in S(n,r)$, and $c\in A(n,r)$ with
$\Delta(c)=\sum_{t}c_{t}\otimes c_{t}'$.

\subsection{The representation $E^{\otimes r}$ of $\Sigma_{r}$}

Let $\{v_{1},\ldots,v_{n}\}$ be a $K$-basis of $E$. Then $E^{\otimes
r}$ has a $K$-basis
$\{v_{\underline{i}}=v_{i_{1}}\otimes\ldots\otimes
v_{i_{r}}|\underline{i}=(i_{1},\ldots,i_{r})\in I(n,r)\}$. The
symmetric group $\Sigma_{r}$ acts on $E^{\otimes r}$ on the right by
$v_{\underline{i}}\sigma=v_{\underline{i}\sigma}$. Then $S(n,r)\cong
\text{End}_{K\Sigma_{r}}(E^{\otimes r})$.

Moreover, if $B$ is the  quotient of $K\Sigma_{r}$ acting faithfully
on $E^{\otimes r}$, then $\text{End}_{S(n,r)}(E^{\otimes r})\cong
B$. When $n\geq r$, $B=K\Sigma_{r}$. This is called Schur--Weyl
duality. (For details, see~\cite{SK}.)

\subsection{The Lie algebra $\mathfrak{gl}_{n}$}\label{S:gln}

Let $K=\mathbb{C}$ and
$\mathfrak{gl}_{n}=\mathfrak{gl}_{n}(\mathbb{C})$ the complex
general linear Lie algebra. Then $\mathcal{U}(\mathfrak{gl}_{n})$,
the universal enveloping algebra of $\mathfrak{gl}_{n}$, acts
naturally on $E$, and on $E^{\otimes r}$ via the comultiplication.
The image of the corresponding representation map is isomorphic to
$S(n,r)$. Therefore there is a surjective homomorphism $\pi$ from
$\mathcal{U}(\mathfrak{gl}_{n})$ to $S(n,r)$. Moreover, the
restriction of this homomorphism to
$\mathcal{U}(\mathfrak{sl}_{n})$, the universal enveloping algebra
of the complex special linear Lie algebra $\mathfrak{sl}_{n}$, is
surjective as well and compatible with $\text{det}^{*}$.

\section{The semigroups $\widetilde{GL}_{n,a}$ and
$\widetilde{SL}_{n,a}$}\label{S:semigroups}

In this section we will introduce two semigroups which will be used
to define the affine Schur algebra.

 Fix $n\in\mathbb{N}$. Let
$\mathfrak{M}_{n}=\mathfrak{M}_{n}(K)=\{M=(m_{ij})_{i,j\in
\mathbb{Z}}|m_{ij}=m_{i+n,j+n}\in K$, and there are only finitely
many nonzero entries in each row of $M \}$. Then with respect to
matrix multiplication $\mathfrak{M}_{n}$ is an algebra with identity
element $I$ whose diagonal entries are $1$ and off-diagonal entries
are $0$. Note that for $M\in\mathfrak{M}_{n}$ there are only
finitely many nonzero entries in each column of $M$.

For $i$, $j\in\mathbb{Z}$, let $E_{ij}\in \mathfrak{M}_{n}$ be the
matrix whose $(i+ln,j+ln)$ entry is $1$ for all $l\in \mathbb{Z}$
and other entries are $0$. Then $E_{ij}=E_{i+n,j+n}$, and
$\{E_{ij}|1\leq i\leq n,j\in\mathbb{Z}\}$ is a $K$-basis of
$\mathfrak{M}_{n}$. The subalgebra of $\mathfrak{M}_{n}$ with basis
$\{E_{ij}|1\leq i,j\leq n\}$ is canonically isomorphic to
$\mathcal{M}_{n}$. We will identify these two algebras. Moreover,
$\mathfrak{M}_{n}\cong\mathcal{M}_{n}\otimes K[t,t^{-1}]$ as
$K$-algebras, and we will identify these two algebras as well. The
isomorphism is given by $E_{i,j+ln}\mapsto E_{ij}\otimes t^{l}$ for
$i,j=1,\ldots,n$, and $l\in\mathbb{Z}$:
\[
\begin{array}{rcl}
E_{i,j+ln}E_{p,q+kn}&=&\delta_{jp}E_{i,q+ln+kn}\\
(E_{ij}\otimes t^{l})(E_{pq}\otimes
t^{k})&=&\delta_{jp}E_{iq}\otimes t^{l+k}.
\end{array}
\]

For $a\in K^{\times}$, $s\in\mathbb{Z}$, tensoring the identity map
of $\mathcal{M}_{n}$ with the $K$-algebra endomorphism $t\mapsto
at^{s}$ of $K[t,t^{-1}]$ defines a $K$-algebra endomorphism
$\eta_{a,s}$ of $\mathfrak{M}_{n}$. Precisely,
$\eta_{a,s}(E_{i,j+ln})=a^{l}E_{i,j+sln}$ for $i,j=1,\ldots,n$, and
$l\in\mathbb{Z}$.

\begin{Lem}\label{L:eta} Let $a,a'\in K^{\times}$, $s,s'\in\mathbb{Z}$.

{\rm (i)} We have $\eta_{a,s}\circ\eta_{a',s'}=\eta_{a'a^{s'},ss'}$.

{\rm (ii)} For $g\in\mathfrak{M}_{n}$, we have
$\eta_{a,s}(g)^{tr}=\eta_{a^{-1},s}(g^{tr})$ where $g^{tr}$ is the
transpose of $g$.

{\rm (iii)} The map $\eta_{a,s}$ fixes elements in
$\mathcal{M}_{n}$.

{\rm (iv)} If $s\not=0$, then $\eta_{a,s}$ is injective. Moreover,
$\eta_{a,\pm 1}$ are isomorphisms.

{\rm (v)} We have $\eta_{a,0}(\mathfrak{M}_{n})=\mathcal{M}_{n}$. We
denote by $\eta_{a}$ this map from $\mathfrak{M}_{n}$ to
$\mathcal{M}_{n}$.
\end{Lem}
\begin{pf}
(i) For $i,j=1,\ldots,n$, $l\in\mathbb{Z}$,
\[
\begin{array}{rcl}
\eta_{a,s}\circ\eta_{a',s'}(E_{i,j+ln})\hspace{-7pt}&=\hspace{-7pt}&\eta_{a,s}((a')^{l}E_{i,j+s'ln})
=(a')^{l}a^{s'l}E_{i,j+ss'ln}\\
\hspace{-7pt}&=\hspace{-7pt}&\eta_{a'a^{s'}}(E_{i,j+ln}).
\end{array}
\]

(ii) For $i,j=1,\ldots,n$, $l\in\mathbb{Z}$,
\[\begin{array}{rcl}
\eta_{a,s}(E_{i,j+ln})^{tr} \hspace{-7pt}& =\hspace{-7pt} &
(a^{l}E_{i,j+sln})^{tr}=a^{l}E_{j,i-sln}\\
\hspace{-7pt}& = \hspace{-7pt}&
\eta_{a^{-1},s}(E_{j,i-ln})=\eta_{a^{-1},s}(E_{i,j+ln}^{tr}).
\end{array}
\]

(iii), (iv) and (v) are clear. \hfill$\square$
\end{pf}

For the rest of this section, we fix an $a\in K^{\times}$. Let
$\widetilde{\text{det}}_{a}=\text{det}\circ \eta_{a}:
\mathfrak{M}_{n}\rightarrow K$. The function
$\widetilde{\text{det}}_{a}$ is a multiplicative function.

Let $\widetilde{GL}_{n,a}=\{M\in\mathfrak{M}_{n} |
\widetilde{\text{det}}_{a}M\not= 0\}$ and
$\widetilde{SL}_{n,a}=\{M\in\mathfrak{M}_{n} |
\widetilde{\text{det}}_{a}M=1\}$. These are two semigroups
containing $GL_{n}$ and $SL_{n}$ respectively.

\begin{Lem}\label{L:etarestricted}
{\rm (i)} For any $b\in K^{\times}$, $b\not=a$, there exists
$g\in\widetilde{GL}_{n,a}$ such that
$\widetilde{\rm{det}}_{b}(g)=0$.

{\rm (ii)} Let $s\in\mathbb{Z}$. The restriction map
$\eta_{a,s}=\eta_{a,s}|_{\widetilde{GL}_{n,a}}:\widetilde{GL}_{n,a}\rightarrow
\widetilde{GL}_{n,1}$ is a semigroup homomorphism fixing elements in
$GL_{n}$. Moreover, it is injective if $s\neq 0$ and bijective if
$s=\pm 1$. The restriction map
$\eta_{a}=\eta_{a}|_{\widetilde{GL}_{n,a}}:
\widetilde{GL}_{n,a}\rightarrow GL_{n}$ is a surjective semigroup
homomorphism fixing elements in $GL_{n}$.

\end{Lem}

\begin{pf}
(i) Let $g=\sum_{i=1}^{n}(E_{ii}-bE_{i,i-n})$. Then
$\widetilde{\text{det}}_{a}(g)=(1-\frac{b}{a})^{n}\not= 0$, but
$\widetilde{\text{det}}_{b}(g)=0$.

(ii) For $g\in\widetilde{GL}_{n,a}$, it follows from
Lemma~\ref{L:eta} (i) that
\[\widetilde{\text{det}}_{1}(\eta_{a,s}(g))=\text{det}\circ\eta_{1}\circ\eta_{a,s}(g)=\text{det}\circ\eta_{a}(g)=\widetilde{\text{det}}_{a}(g)\not=0.\]
Therefore $\eta_{a,s}(g)\in\widetilde{GL}_{n,1}$. The other
statements follow directly from Lemma~\ref{L:eta}.
 \hfill$\square$
\end{pf}
\vskip5pt

Let $G_{1}$ be the subgroup of $\widetilde{GL}_{n,a}$ generated by
$s_{i}=E_{i,i+1}+E_{i+1,i}+\sum_{j\not= i,i+1}E_{jj}$ ($\sum_{j}$
means $\sum_{j=1}^{n}$) where $i=1,\ldots,n-1$. Then $G_{1}$ is
contained in $GL_{n}$ and $G_{1}$ is isomorphic to the symmetric
group $\Sigma_{n}$. Let $G_{2}$ be the subgroup of
$\widetilde{GL}_{n,a}$ generated by $\tau_{i}=E_{i,i+n}+\sum_{j\not=
i}E_{jj}$ where $i=1,\ldots,n$. Note that an element $\tau\in
\widetilde{GL}_{n,a}$ is in $G_{2}$ if and only if $\tau$ is of the
form $\tau=\sum_{s=1}^{n}E_{s,s+n\varepsilon_{s}}$ where
$\varepsilon_{1},\ldots,\varepsilon_{n}\in\mathbb{Z}$. Therefore
$G_{2}$ and $G_{1}$ intersect trivially since no elements in $G_{2}$
lies in $GL_{n}$ except the identity matrix. Moreover
$\sum_{s=1}^{n}E_{s,s+n\varepsilon_{s}}\mapsto
(\varepsilon_{1},\ldots,\varepsilon_{n})$ is an isomorphism from
$G_{2}$ to $\mathbb{Z}^{n}$. The group $G_{1}$ acts on $G_{2}$ on
the right by conjugation
\[\begin{array}{rcl}s_{i}(\sum_{s=1}^{n}E_{s,s+n\varepsilon_{s}})s_{i}&=&E_{i,i+n\varepsilon_{i+1}}+E_{i+1,i+1+n\varepsilon_{i}}+\sum_{s\neq
i,i+1}E_{s,s+n\varepsilon_{s}}\\
& = & \sum_{s=1}^{n}E_{s,s+n\varepsilon_{s_{i}(s)}}.\end{array}\]
Taking the isomorphisms given above into account this right action
coincides with the right action of $\Sigma_{n}$ on $\mathbb{Z}^{n}$
by place permutation. Therefore we have

\begin{Lem}\label{L:equivgp}
The subgroup $\langle G_{1},G_{2}\rangle$ of $\widetilde{GL}_{n,a}$
generated by $G_{1}$ and $G_{2}$ is the semi-direct product
 $G_{1}\ltimes G_{2}$ of $G_{1}$ and $G_{2}$, and hence isomorphic to
$\widehat{\Sigma}_{n}=\Sigma_{r}\ltimes\mathbb{Z}^{n}$, the extended
affine Weyl group (of type $A$). \hfill$\square$
\end{Lem}

We will identify $G_{1}$ with $\Sigma_{n}$, $G_{2}$ with
$\mathbb{Z}^{n}$ and $\langle G_{1},G_{2}\rangle$ with
$\widehat{\Sigma}_{n}$.  Consequently, the group
$\widehat{\Sigma}_{n}$ acts on $\mathfrak{M}_{n}$ by conjugation.
The group $\widehat{\Sigma}_{n}$ has another presentation
(see~\cite{RMG}):

generators : $s_{1},\ldots,s_{n-1},s_{n},\rho$, where
$s_{n}=E_{10}+E_{n,n+1}+\sum_{j\not= 1,n}E_{jj}$,
$\rho=\sum_{j=1}^{n}E_{j,j+1}$,

relations :
\[\begin{array}{rcl}
s_{i}^2\hspace{-7pt} & =\hspace{-7pt} & I,~~\text{for}~i=1,\ldots,n\\
s_{i}s_{\overline{i+1}}s_{i} \hspace{-7pt}& = \hspace{-7pt}& s_{\overline{i+1}}s_{i}s_{\overline{i+1}},~~\text{for}~i=1,\ldots,n\\
s_{j}s_{i} \hspace{-7pt}& =\hspace{-7pt} &
s_{i}s_{j},~~\text{for}~i,j=1,\ldots,n,~j\not= \overline{i\pm
1}\\
\rho s_{\overline{i+1}}\rho^{-1}\hspace{-7pt}& =\hspace{-7pt} &
s_{i},~~\text{for}~i=1,\ldots,n
\end{array}\]
where \ $\bar{}: \mathbb{Z}\rightarrow \{1,\ldots,n\}$ is the map
taking least positive remainder modulo $n$. We claim that for
$w\in~\widehat{\Sigma}_{n}$ and $g=(g_{ij})\in\mathfrak{M}_{n}$
\[wgw^{-1}=(g_{ij}'=g_{w^{-1}(i),w^{-1}(j)})\]
where for $z\in\mathbb{Z}$ and $i=1,\ldots,n$
\[
\begin{array}{rcl}
\rho(z)&=&z-1,\\
s_{i}(z)&=&\begin{cases}z &   {\rm \ if\ } z\not\equiv
i,i+1\pmod{n},\\
z+1 & {\rm\ if\ } z\equiv i\pmod{n},\\
z-1 & {\rm\ if\ } z\equiv i+1\pmod{n}.\end{cases} \end{array}\] The
proof follows by checking on generators.

In particular,

\begin{Lem}\label{L:weylaction}
The extended affine Weyl group $\widehat{\Sigma}_{n}$ acts on
$\widetilde{GL}_{n,a}$ and $\widetilde{SL}_{n,a}$:
$w\in\widehat{\Sigma}_{n}$ sends $g=(g_{ij})$ to
$g'=wgw^{-1}=(g'_{ij})$ where $g'_{ij}=g_{w^{-1}(i),w^{-1}(j)}$.
\hfill$\square$
\end{Lem}
\vskip10pt

Let $K^{\widetilde{GL}_{n,a}}$ be the set of maps from
$\widetilde{GL}_{n,a}$ to $K$. It is a commutative algebra with
pointwise multiplication $(ff')(g)=f(g)f'(g)$ for $f,f'\in
K^{\widetilde{GL}_{n,a}}$ and $g\in\widetilde{GL}_{n,a}$. For
$i,j\in \mathbb{Z}$, let $c_{ij}$ be the coordinate function
$c_{ij}:\widetilde{GL}_{n,a}\rightarrow K$, $g=(g_{pq})\mapsto
g_{ij}$. Then
$c_{ij}=c_{i+n,j+n}=c_{\overline{i},j+\overline{i}-i}$.

Let $I(\mathbb{Z},r)=\{\underline{i}=(i_{1},\ldots,i_{r})|i_{t}\in
\mathbb{Z}, t=1,\ldots,r\}$. Then the extended affine Weyl group
$\widehat{\Sigma}_{r}$ acts on $I(\mathbb{Z},r)$ on the right with
$\Sigma_{r}$ acting by place permutation and $\mathbb{Z}^{r}$ acting
by shifting, i.e.
$\underline{i}(\sigma,\varepsilon)=\underline{i}\sigma+n\varepsilon$
for $\underline{i}\in I(\mathbb{Z},r),~\sigma\in\Sigma_{r}$ and
$\varepsilon\in\mathbb{Z}^{r}$. The group $\widehat{\Sigma}_{r}$
acts on $I(\mathbb{Z},r)\times I(\mathbb{Z},r)$ diagonally. The
equivalence relations on $I(\mathbb{Z},r)$ and on
$I(\mathbb{Z},r)\times I(\mathbb{Z},r)$ induced by this action will
be both denoted by the symbol $\sim_{\widehat{\Sigma}_{r}}$. For
$\underline{i}=(i_{1},\ldots,i_{r}),\underline{j}=(j_{1},\ldots,j_{r})\in
I(\mathbb{Z},r)$ write
$c_{\underline{i},\underline{j}}=c_{i_{1}j_{1}}\cdots
c_{i_{r}j_{r}}$, then
$c_{\underline{i},\underline{j}}=c_{\underline{p},\underline{q}}$ if
$(\underline{i},\underline{j})\sim_{\widehat{\Sigma}_{r}}
(\underline{p},\underline{q})$.

Let $K[T]=K[t_{ij}]_{i=1,\ldots,n,j\in\mathbb{Z}}$ be the polynomial
algebra in indeterminates $t_{ij},i=1,\ldots,n,j\in\mathbb{Z}$. For
an integer $r\geq 0$ let $K[T]_{r}$ denote the homogeneous component
of $K[T]$ of degree $r$, and let $\overline{K[T]_{r}}$ denote the
closure of $K[T]_{r}$ with respect to the basis
$t_{\underline{i},\underline{j}}=t_{i_{1}j_{1}}\cdots
t_{i_{r}j_{r}}, \underline{i}\in I(n,r),\underline{j}\in
I(\mathbb{Z},r)$. We have
$t_{\underline{i},\underline{j}}=t_{\underline{p},\underline{q}}$ if
and only if
$(\underline{i},\underline{j})\sim_{\Sigma_{r}}(\underline{p},\underline{q})$.

\begin{Prop}\label{P:algindep3} If $P$ is a nonzero element in $\bigoplus_{r\geq
0}\overline{K[T]_{r}}$, then there exists $g\in\widetilde{GL}_{n,a}$
such that $P(c_{ij})(g)\neq 0$.\hfill$\square$
\end{Prop}

\begin{pf} For $r\geq 0$ let $I_{r}$ be a fixed set of representatives of $(I(n,r)\times
I(n,r))/\Sigma_{r}$, and for $(\underline{i},\underline{j})\in
I_{r}$ let $Z_{\underline{i},\underline{j}}$ be a fixed set of
representatives of
$\mathbb{Z}^{r}/\Sigma_{\underline{i},\underline{j}}$. Then
$\{(\underline{i},\underline{j}+n\varepsilon)|(\underline{i},\underline{j})\in
I_{r},\varepsilon\in Z_{\underline{i},\underline{j}}\}$ is a set of
representatives of $(I(n,r)\times I(\mathbb{Z},r))/\Sigma_{r}$.
Suppose
\[P=\sum_{r\geq 0}\sum_{(\underline{i},\underline{j})\in I_{r}}\sum_{\varepsilon\in Z_{\underline{i},\underline{j}}}
\lambda_{\underline{i},\underline{j},\varepsilon}
t_{\underline{i},\underline{j}+n\varepsilon}.\]

Since $P$ is not zero, there exists $r_{0}$,
$(\underline{i}^{0},\underline{j}^{0})\in I_{r_{0}}$ and
$\varepsilon^{0}\in Z_{\underline{i}^{0},\underline{j}^{0}}$ such
that
$\lambda_{\underline{i}^{0},\underline{j}^{0},\varepsilon^{0}}\neq
0$. Set
$L=\{l_{1},\ldots,l_{s}\}=\{\varepsilon^{0}_{1},\ldots,\varepsilon^{0}_{r_{0}}\}$,
where $l_{1},\ldots,l_{s}$ are distinct.

Let $X_{l_{1}},\ldots,X_{l_{s}}$ be a set of indeterminates labeled
by the set $L$. Consider
\begin{equation}\label{E:algindep3}
\begin{cases}a^{l_{1}}X_{l_{1}}+\ldots+a^{l_{s}}X_{l_{s}}-1=0\\
\sum_{\varepsilon\in I(L,r_{0})\cap
Z_{\underline{i}^{0},\underline{j}^{0}}}\lambda_{\underline{i}^{0},\underline{j}^{0},\varepsilon}X_{\varepsilon_{1}}\cdots
X_{\varepsilon_{r_{0}}}\neq 0.\end{cases}
\end{equation}
Since $K$ is infinite, the subvariety of $K^{n}$ defined by the
equations (\ref{E:algindep3}) is nonempty (by the first equation we
can express $X_{l_{s}}$ as the sum of a linear combination of
$X_{l_{1}},\ldots,X_{l_{s-1}}$ and a nonzero constant. Substituting
it into the second equation we see that the resulting polynomial is
not zero. Then induct on $s$). That is, there exist elements
$x_{l_{1}},\ldots,x_{l_{s}}$ in $K$ satisfying the equations
(\ref{E:algindep3}). Set
\[P'=\sum_{r\geq 0}\sum_{(\underline{i},\underline{j})\in I_{r}}\big(\sum_{\varepsilon\in I(L,r)\cap Z_{\underline{i},\underline{j}}}
\lambda_{\underline{i},\underline{j},\varepsilon}x_{\varepsilon_{1}}\cdots
x_{\varepsilon_{r}}\big) t_{\underline{i},\underline{j}}.\] Then
$P'$ is a nonzero polynomial in indeterminates $t_{ij}$
($i,j=1,\ldots,n$). Therefore there exists $g'\in GL_{n}$ such that
$P'(c_{ij})(g')\neq 0$.

Define $g=(g_{pq})$ to be the matrix in $\mathfrak{M}_{n}$ whose
$(i,j+ln)$ ($i,j=1,\ldots,n,l\in\mathbb{Z}$) entry is $x_{l}g_{ij}'$
if $l\in L$ and $0$ otherwise.
\[{\setlength{\unitlength}{0.7pt}
\begin{picture}(450,60)
\drawline(0,0)(430,0) \drawline(0,60)(430,60) \drawline(50,0)(50,60)
\drawline(110,0)(110,60)
\put(65,25){$x_{l_{1}}g'$}\put(75,-15){\footnotesize$l_{1}$}
\drawline(170,0)(170,60)\put(135,25){$0$} \drawline(230,0)(230,60)
\put(185,25){$x_{l_{2}}g'$}\put(195,-15){\footnotesize$l_{2}$}
\drawline(290,0)(290,60)\put(255,25){$0$}\drawline(350,0)(350,60)\put(315,25){$0$}
\drawline(410,0)(410,60)\put(365,25){$x_{l_{3}}g'$}\put(375,-15){\footnotesize$l_{3}$}
\end{picture}}\]
Therefore
\[c_{\underline{i},\underline{j}+n\varepsilon}(g)
=\begin{cases}x_{\varepsilon_{1}}\cdots x_{\varepsilon_{r}}c_{\underline{i},\underline{j}}(g'), & \text{if}~\varepsilon\in I(L,r),\\
                             0, & \text{otherwise}.
                             \end{cases}\]
Therefore
\[\begin{array}{rcl}
P(c_{ij})(g)\hspace{-7pt} &=\hspace{-7pt}&\sum_{r\geq
0}\sum_{(\underline{i},\underline{j})\in I_{r}}\sum_{\varepsilon\in
Z_{\underline{i},\underline{j}}}
\lambda_{\underline{i},\underline{j},\varepsilon}
c_{\underline{i},\underline{j}+n\varepsilon}(g)\\
\hspace{-7pt}&=\hspace{-7pt}&\sum_{r\geq
0}\sum_{(\underline{i},\underline{j})\in I_{r}}\sum_{\varepsilon\in
I(L,r)\cap Z_{\underline{i},\underline{j}}}
\lambda_{\underline{i},\underline{j},\varepsilon}
c_{\underline{i},\underline{j}+n\varepsilon}(g)\\
\hspace{-7pt}&=\hspace{-7pt}&\sum_{r\geq
0}\sum_{(\underline{i},\underline{j})\in
I_{r}}\big(\sum_{\varepsilon\in I(L,r)\cap
Z_{\underline{i},\underline{j}}}
\lambda_{\underline{i},\underline{j},\varepsilon}x_{\varepsilon_{1}}\cdots
x_{\varepsilon_{r}}\big) c_{\underline{i},\underline{j}}(g)\\
\hspace{-7pt}&=\hspace{-7pt}& P'(c_{ij})(g')\neq 0.\end{array}\]

By definition we have $g=\sum_{i,j=1,\ldots,n,l\in
L}x_{l}g_{ij}'E_{i,j+ln}$. Therefore
\[\eta_{a}(g)=\sum_{i,j=1,\ldots,n}(\sum_{l\in
L}a^{l}x_{l})g_{ij}'E_{ij}=g',\] which implies
$g\in\widetilde{GL}_{n,a}$. \hfill$\square$
\end{pf}

If $P$ in Proposition~\ref{P:algindep3} is homogeneous, i.e. $P\in
\overline{K[T]_{r}}$ for some $r$, then the polynomial $P'$ in the
proof is homogeneous of degree $r$. In this case, we can choose $g'$
to be in $SL_{n}$, and then the resulting $g$ lies in
$\widetilde{SL}_{n,a}$. Namely, we have the following proposition.

\begin{Prop}\label{P:algindep2} Let $r\geq 0$ be an integer. If $P$ is a nonzero element in
$\overline{K[T]_{r}}$, then there exists $g\in\widetilde{SL}_{n,a}$
such that $P(c_{ij})(g)\neq 0$. \hfill$\square$
\end{Prop}

By Proposition~\ref{P:algindep3} the functions
$\{c_{ij}|i=1,\ldots,n,j\in\mathbb{Z}\}$ are algebraically
independent, and hence the subalgebra $\widetilde{A}(n)$ of
$K^{\widetilde{GL}_{n,a}}$ generated by all $c_{ij}$'s is a
polynomial algebra in indeterminates $\{c_{ij}\}$, where $(i,j)$
runs over a set of representatives of $\widehat{\Sigma}_{1}$-orbits
of $\mathbb{Z}\times\mathbb{Z}$ (here
$\widehat{\Sigma}_{1}=\mathbb{Z}$ acts on
$\mathbb{Z}=I(\mathbb{Z},1)$ by shifting and on
$\mathbb{Z}\times\mathbb{Z}$ diagonally), for example we may take
$\{(i,j)|i=1,\ldots,n,j\in\mathbb{Z}\}$. Let $r\geq 0$ be an
integer. Denote by $\widetilde{A}(n,r)$ the homogeneous component of
$\widetilde{A}(n)$ of degree $r$. Note that
$c_{\underline{i},\underline{j}}=c_{\underline{p},\underline{q}}$ if
and only if
$(\underline{i},\underline{j})\sim_{\widehat{\Sigma}_{r}}(\underline{p},\underline{q})$.
Thus $\widetilde{A}(n,r)$ has a basis
$\{c_{\underline{i},\underline{j}}\}$, where
$(\underline{i},\underline{j})$ runs over a set of representatives
of $\widehat{\Sigma}_{r}$-orbits of $I(\mathbb{Z},r)\times
I(\mathbb{Z},r)$. By Proposition~\ref{P:algindep3} the space
$\bigoplus_{r\geq 0}\overline{\widetilde{A}(n,r)}$ is a subspace of
$K^{\widetilde{GL}_{n,a}}$.

\section{The affine Schur algebra $\widetilde
S(n,r)$}\label{S:affineschuralgebra}

 Fix $n$, $r\in {\mathbb{N}}$. In this section we will define the
 affine Schur algebra and show how it is related to the two
 semigroups $\widetilde{GL}_{n,a}$ and $\widetilde{SL}_{n,a}$ and
 the extended affine Weyl group $\widehat{\Sigma}_{r}$.

Let, for a moment, $G$ be any semigroup with identity $1_{G}$. Let
$K^{G}$ be the $K$-space of all maps from $G$ to $K$. Then $K^{G}$
is a commutative $K$-algebra. The semigroup structure on $G$ gives
rise to two maps
\[
\begin{array}{c}\Delta=\Delta_{G}:K^{G}\rightarrow
K^{G\times G},~~ f\mapsto
((s,t)\mapsto f(st) {\rm\ for\ } s,t\in G )\\
\epsilon=\epsilon_{_G}:K^{G}\rightarrow K,~~ f\mapsto
f(1_{G}).\end{array}\] Both $\Delta$ and $\epsilon$ are $K$-algebra
homomorphisms (see~\cite{JAG}).

The restrictions of $\Delta$ and $\epsilon$ to $\widetilde{A}(n,r)$
can be written explicitly as follows
\[\begin{array}{c}\Delta: \widetilde A(n,r)\rightarrow \overline{\widetilde
A(n,r) \otimes \widetilde A(n,r)},~~
c_{\underline{i},\underline{j}}\mapsto \sum_{\underline{s}\in
I(\mathbb{Z},r)}{c_{\underline{i},\underline{s}}\otimes
c_{\underline{s},\underline{j}} }\\
\epsilon:\widetilde{A}(n,r)\rightarrow K,~~
c_{\underline{i},\underline{j}}\mapsto
\delta_{\underline{i},\underline{j}}.\end{array}\]

\begin{Lem} These two maps make $\widetilde{A}(n,r)$
into a formal $K$-coalgebra (for definition of formal coalgebras see
Appendix 1).
\end{Lem}
\begin{pf} One can check by the definition of formal
coalgebras. \hfill$\square$
\end{pf}

It follows by Theorem~\ref{T:fromcoalgebratoalgebra} that
$\widetilde S(n,r)=\widetilde A(n,r)^{\#}$ is a $K$-algebra, called
{\em the affine Schur algebra}. Let
$\{\xi_{\underline{i},\underline{j}}|\underline{i},\underline{j}\in
I(\mathbb{Z},r)\}$ be the basis dual to
$\{c_{\underline{i},\underline{j}}|\underline{i},\underline{j}\in
I(\mathbb{Z},r)\}$, i.e.
\[
\xi_{\underline{i},\underline{j}}(c_{\underline{p},\underline{q}})=\begin{cases}
1, & {\rm
if}\ \ (\underline{i},\underline{j})\sim_{\widehat{\Sigma}_{r}}(\underline{p},\underline{q})\\
0, & {\rm otherwise}.\end{cases}
\]
Therefore
$\xi_{\underline{i},\underline{j}}=\xi_{\underline{p},\underline{q}}$
if and only if
$(\underline{i},\underline{j})\sim_{\widehat{\Sigma}_{r}}
(\underline{p},\underline{q})$.

Recall that \ $\bar{}: \mathbb{Z}\rightarrow \{1,\ldots,n\}$ is the
map taking least positive remainder modulo $n$. It can be extended
to \ $\bar{}: I(\mathbb{Z},r)\rightarrow I(n,r)$. Then
$\xi_{\underline{i},\underline{j}}=\xi_{\overline{\underline{i}},\underline{j}+\overline{\underline{i}}-\underline{i}}=\xi_{\underline{i}+\overline{\underline{j}}-\underline{j},\overline{\underline{j}}}$\
. In some later cases we will assume that $\underline{i}\in I(n,r)$
or $\underline{j}\in I(n,r)$.

For $\underline{i},\underline{j},\underline{k},\underline{l}\in
I(\mathbb{Z},r)$, we have the following formula, known as Schur's
product rule:
\[
\xi_{\underline{i},\underline{j}}\xi_{\underline{k},\underline{l}}=\sum_{(\underline{p},\underline{q})\in
( I(\mathbb{Z},r)\times
I(\mathbb{Z},r))/\widehat{\Sigma}_{r}}{Z(\underline{i},\underline{j},\underline{k},\underline{l},\underline{p},\underline{q})\xi_{\underline{p},\underline{q}}}
\]
where
$Z(\underline{i},\underline{j},\underline{k},\underline{l},\underline{p},\underline{q})=
\#\{\underline{s}\in
I(\mathbb{Z},r)|(\underline{i},\underline{j})\sim_{\widehat{\Sigma}_{r}}(\underline{p},\underline{s}),(\underline{s},\underline{q})\sim_{\widehat{\Sigma}_{r}}(\underline{k},\underline{l})\}$.

Directly from Schur's product rule, we have the following
proposition.

\begin{Prop}
{\rm(i)} We have
 $\xi_{\underline{i},\underline{j}}\xi_{\underline{k},\underline{l}}=0$
unless $\underline{j} \sim_{\widehat{\Sigma}_{r}} \underline{k}$.

{\rm(ii)} We have
$\xi_{\underline{i},\underline{i}}\xi_{\underline{i},\underline{j}}=\xi_{\underline{i},\underline{j}}=\xi_{\underline{i},\underline{j}}\xi_{\underline{j},\underline{j}}$,
for $\underline{i},\underline{j}\in I(\mathbb{Z},r)$.

{\rm(iii)} $\sum_{\underline{i}\in
I(n,r)/\Sigma_{r}}{\xi_{\underline{i},\underline{i}}}$ is a
decomposition of unity into orthogonal idempotents.

{\rm(iv)} The subalgebra of $\widetilde{S}(n,r)$ with basis
$\{\xi_{\underline{i},\underline{j}}|\underline{i},\underline{j}\in
I(n,r)\}$ is naturally isomorphic to $S(n,r)$. We will identify
these two algebras.

\end{Prop}

Let $\widetilde{e}_{a}: K\widetilde{GL}_{n,a}\rightarrow
\widetilde{S}(n,r)$ be the algebra homomorphism sending
$g\in\widetilde{GL}_{n,a}$ to $(\widetilde{e}_{a}(g): c\mapsto
c(g))$ for any $c\in \widetilde{A}(n,r)$. The image of $KGL_{n}$
under $\widetilde{e}_{a}$ lies in $S(n,r)$. In fact, the restriction
map $\widetilde{e}_{a}|_{KGL_{n}}$ is the map $e$ defined in
Section~\ref{S:naturalrep}. The following theorem follows from
Proposition~\ref{P:algindep2}.

\begin{Thm}\label{T:evnsurj} The map
$\widetilde{e}_{a}$ is surjective. Moreover, the restriction map
$\widetilde{ e}_{a}|_{K\widetilde{SL}_{n,a}}:
K\widetilde{SL}_{n,a}\rightarrow \widetilde{S}(n,r)$ is surjective.
\end{Thm}
\begin{pf} Suppose the image of the restriction map $\widetilde{
e}_{a}|_{K\widetilde{SL}_{n,a}}$ is a proper subspace of
$\widetilde{S}(n,r)$. Then there exists a nonzero
$c\in\overline{\widetilde{A}(n,r)}=S(n,r)^{*}$ such that $c(\xi)=0$
for any $\xi\in
\text{Im}~\widetilde{e}_{a}|_{K\widetilde{SL}_{n,a}}$. In
particular, $c(g)=c(\widetilde{e}_{a}(g))=0$ for any
$g\in\widetilde{SL}_{n,a}$, contradicting
Proposition~\ref{P:algindep2}. \hfill$\square$
\end{pf}

Let $\widetilde{Y}_{a}$ be the kernel of $\widetilde{e}_{a}$, i.e.
$0\rightarrow \widetilde{Y}_{a} \rightarrow
K\widetilde{GL}_{n,a}\stackrel{\widetilde{e}_{a}}{\rightarrow}
\widetilde{S}(n,r)\rightarrow 0$ is exact.

\begin{Prop}\label{P:cf}
 Let $f\in K^{\widetilde{GL}_{n,a}}$. Then
$f\in\overline{\widetilde{A}(n,r)}$ if and only if \
$f(\widetilde{Y}_{a})=0$.
\end{Prop}
\begin{pf} Assume
$f\in\overline{\widetilde{A}(n,r)}$. For any $y\in
\widetilde{Y}_{a}$, we have $\widetilde{e}_{a}(y)=0$, i.e.
$\widetilde{e}_{a}(y)(c_{\underline{i},\underline{j}})=c_{\underline{i},\underline{j}}(y)=0$
for any $\underline{i}\in I(n,r)$, $\underline{j}\in
I(\mathbb{Z},r)$. Hence $f(y)=0$.

Assume $f(\widetilde{Y}_{a})=0$. Then there exists $c\in
\overline{\widetilde{A}(n,r)}=S(n,r)^{*}$ such that
$c(g)=\widetilde{e}_{a}(g)(c)=c(\widetilde{e}_{a}(g))=f(g)$  for any
$g\in \widetilde{GL}_{n,a}$, and hence $f=c$. \hfill$\square$
\end{pf}

Let $W$ be a representation of $\widetilde{GL}_{n,a}$ with basis
$\{w_{j}|j\in J\}$. Suppose $g(w_{j})=\sum_{j'\in
J}r_{j'j}(g)w_{j'}$, for $g\in\widetilde{GL}_{n,a}$, $j\in J$. We
say that $W$ is a representation with coefficients in
$\overline{\widetilde{A}(n,r)}$ if $r_{j'j}\in
\overline{\widetilde{A}(n,r)}$ for any $j,j'\in J$. The matrix
$(r_{j'j})_{j,j'\in J}$ is called the coefficient matrix of $W$. Let
$\widetilde{M}(n,r)$ denote the category of representations of
$\widetilde{GL}_{n,a}$ with coefficients in
$\overline{\widetilde{A}(n,r)}$.

\begin{Prop} \label{P:repreequiv}
 Let $W$ be a representation of
$\widetilde{GL}_{n,a}$. Then $W\in \widetilde{M}(n,r)$ if and only
if $\widetilde{Y}_{a}W=0$. This induces an equivalence between
$\widetilde{M}(n,r)$ and $\widetilde{S}(n,r)$-$\rm{Mod}$.
\end{Prop}

\begin{pf} Let $\{w_{j}|j\in J\}$ be the basis
of $W$, and $(r_{j'j})_{j,j'\in J}$ the coefficient matrix of $W$
with respect to this basis. Then $W\in \widetilde{M}(n,r)$ if and
only if $r_{j'j}\in \overline{\widetilde{A}(n,r)}$. By
Proposition~\ref{P:cf}, this is equivalent to
$\widetilde{Y}_{a}W=0$. \hfill$\square$
\end{pf}
\vskip10pt Let $\widetilde{E}=K\{v_{i}|i\in \mathbb{Z}\}$. Then
$\widetilde{E}^{\otimes
r}=K\{v_{\underline{i}}=v_{i_{1}}\otimes\ldots\otimes v_{i_{r}}|
\underline{i}=(i_{1},\ldots,i_{r})\in I(\mathbb{Z},r)\}$. The
semigroup $\widetilde{GL}_{n,a}$ acts on $\widetilde{E}$ from the
left by matrix multiplication, and on $\widetilde{E}^{\otimes r}$
diagonally, i.e. $g(v_{\underline{i}})=\sum_{\underline{j}\in
I(\mathbb{Z},r)}{c_{\underline{j},\underline{i}}(g)v_{\underline{j}}}$
for $\underline{i}\in I(\mathbb{Z},r)$. By
Proposition~\ref{P:repreequiv}, the tensor space
$\widetilde{E}^{\otimes r}$ can be regarded as an
$\widetilde{S}(n,r)$-module,
$\xi(v_{\underline{i}})=\sum_{\underline{j}\in
I(\mathbb{Z},r)}{\xi(c_{\underline{j},\underline{i}})v_{\underline{j}}}$.
This is a faithful module. Thus the representation map
$\phi:\widetilde{S}(n,r)\rightarrow
\text{End}_{K}(\widetilde{E}^{\otimes r})$ is injective.

The extended affine Weyl group $\widehat{\Sigma}_{r}$ acts on
$\widetilde{E}^{\otimes r}$ on the right by
$v_{\underline{i}}w=v_{\underline{i}w}$. These two actions commute.
In particular, the image of $\phi$ lies in
$\text{End}_{K\widehat{\Sigma}_{r}}(\widetilde{E}^{\otimes r})$. In
fact we have

\begin{Thm} \label{T:defaffine}
The image of $\phi$ is exactly
$\rm{End}_{K\widehat{\Sigma}_{r}}(\widetilde{E}^{\otimes r})$.
Therefore $\phi$ induces an isomorphism between $\widetilde{S}(n,r)$
and $ \rm{End}_{K\widehat{\Sigma}_{r}}(\widetilde{E}^{\otimes r})$.
\end{Thm}
\begin{pf} For $\underline{i}$, $\underline{j}\in I(\mathbb{Z},r)$, let
$x_{\underline{i},\underline{j}}$ denote the endomorphism of
$\widetilde{E}^{\otimes r}$ sending $v_{\underline{l}}$ to
$\delta_{\underline{j},\underline{l}}v_{\underline{i}}$. Then
\[
\text{End}_{K}(\widetilde{E}^{\otimes
r})=\big\{\sum_{\underline{i},\underline{j}\in
I(\mathbb{Z},r)}{m_{\underline{i},\underline{j}}x_{\underline{i},\underline{j}}}|m_{\underline{i},\underline{j}}\in
K {\rm\ and \ for \ each\ } \underline{j},
m_{\underline{i},\underline{j}}=0 {\rm \ for \ almost \  all\ }
\underline{i}\big\}
\]

Since
\[\xi_{\underline{i},\underline{j}}v_{\underline{l}}=\sum_{\underline{k}\in I(\mathbb{Z},r)}\xi_{\underline{i},\underline{j}}(c_{\underline{k},\underline{l}})v_{\underline{k}}=\sum_{w\in\widehat{\Sigma}_{\underline{i},\underline{j}}\backslash\widehat{\Sigma}_{r}:\underline{j}w=\underline{l}}v_{\underline{i}w},\]
the image of $\xi_{\underline{i},\underline{j}}$ under $\phi$ is
\[\phi(\xi_{\underline{i},\underline{j}})=\sum_{w\in\widehat{\Sigma}_{\underline{i},\underline{j}}\backslash\widehat{\Sigma}_{r}}x_{\underline{i}w,\underline{j}w}.\]

The right action of $\widehat{\Sigma}_{r}$ on
$\widetilde{E}^{\otimes r}$ induces a right action of
$\widehat{\Sigma}_{r}$ on $\text{End}_{K}(\widetilde{E}^{\otimes
r})$,
\[
f^{w}=\sum_{\underline{i},\underline{j}\in
I(\mathbb{Z},r)}{m_{\underline{i},\underline{j}}x_{\underline{i}w,\underline{j}w}}=\sum_{\underline{i},\underline{j}\in
I(\mathbb{Z},r)}{m_{\underline{i}w^{-1},\underline{j}w^{-1}}x_{\underline{i},\underline{j}}}
\]
for $f=\sum_{\underline{i},\underline{j}\in
I(\mathbb{Z},r)}{m_{\underline{i},\underline{j}}x_{\underline{i},\underline{j}}}\in
\text{End}_{K}(\widetilde{E}^{\otimes r}) $ and $w\in
\widehat{\Sigma}_{r}$. It is easy to see $f$ is fixed by
$\widehat{\Sigma}_{r}$ if and only if
$m_{\underline{i},\underline{j}}=m_{\underline{i}w,\underline{j}w}$
for any $\underline{i},\underline{j}\in I(\mathbb{Z},r)$ and $w\in
\widehat{\Sigma}_{r}$, that is, $f$ is a linear combination of
$\phi(\xi_{\underline{i},\underline{j}})$'s. Now the desired result
follows from
$\text{End}_{K\widehat{\Sigma}_{r}}(\widetilde{E}^{\otimes
r})=(\text{End}_{K}(\widetilde{E}^{\otimes
r}))^{\widehat{\Sigma}_{r}}$. \hfill$\square$

\end{pf}
\vskip5pt

We will identify $\phi(\xi_{\underline{i},\underline{j}})$ with
$\xi_{\underline{i},\underline{j}}$. Now let us give an analogue of
J.A.Green's product formula (\ref{E:greenprod})(\cite{JAG3}(2.6))
for the finite Schur algebra.

For $\underline{i}$, $\underline{j}\in I(\mathbb{Z},r)$, let
$\widehat{\Sigma}_{\underline{i}}$ be the stabilizer of
$\underline{i}$ in $\widehat{\Sigma}_{r}$ and
$\widehat{\Sigma}_{\underline{i},\underline{j}}=\widehat{\Sigma}_{\underline{i}}\cap
\widehat{\Sigma}_{\underline{j}}$, and so on. Note that if
$\underline{i}\in I(n,r)$, then
$\widehat{\Sigma}_{\underline{i}}=\Sigma_{\underline{i}}$.

Let
$A=K\{x_{\underline{i},\underline{j}}|\underline{i},\underline{j}\in
I(\mathbb{Z},r)\}$, $B=\text{End}_{K}(\widetilde{E}^{\otimes r})$.
For a subgroup $H$ of $\widehat{\Sigma}_{r}$, let
\[\begin{array}{rcl}A_{H}&=&\{a\in A|a^{h}=a {\rm\ for\ any}\ h\in
H\}\\
B_{H}&=&\{b\in B|b^{h}=b {\rm\ for\ any}\ h\in H\}\end{array}\] then
\[T_{H,\widehat{\Sigma}_{r}}:A_{H}\rightarrow B_{H},~~ a\mapsto\sum_{g\in H\backslash\widehat{\Sigma}_{r}}a^{g}\]
is well-defined (see Appendix 2). For $\delta\in
\widehat{\Sigma}_{r}$ we set $H^{\delta}=\delta^{-1}H\delta$.

The following lemma is a special case of Appendix 2
Lemma~\ref{L:mackey}: Mackey's formula.

\begin{Lem}\label{L:specialmackey}
Let $H_{1}$ and $H_{2}$ be two subgroups of $\widehat{\Sigma}_{r}$,
and $a\in A_{H_{1}}$, $b\in A_{H_{2}}$. If \
$aT_{H_{2},\widehat{\Sigma}_{r}}(b)\in A_{H_{1}}$, $T_{H_{1}\cap
H_{2}^{\delta},H_{1}}(ab^{\delta})=0$ for almost all $\delta$, and
$T_{H_{1}\cap H_{2}^{\delta},H_{1}}(ab^{\delta})\in A_{H_{1}}$ for
all $\delta$, then
\[
T_{H_{1},\widehat{\Sigma}_{r}}(a)T_{H_{2},\widehat{\Sigma}_{r}}(b)=\sum_{\delta\in
H_{2}\backslash\widehat{\Sigma}_{r}/H_{1}}T_{H_{1}\cap
H_{2}^{\delta},\widehat{\Sigma}_{r}}(ab^{\delta}).
\] \hfill$\square$
\end{Lem}

Let $\underline{i},\underline{j},\underline{l}\in
I(\mathbb{Z},r)$. Then
\[\begin{array}{c}
\xi_{\underline{i},\underline{j}}=T_{\widehat{\Sigma}_{\underline{i},\underline{j}},\widehat{\Sigma}_{r}}(x_{\underline{i},\underline{j}})\ \\
\xi_{\underline{j},\underline{l}}=T_{\widehat{\Sigma}_{\underline{j},\underline{l}},\widehat{\Sigma}_{r}}(x_{\underline{j},\underline{l}}).
\end{array}\]

Therefore setting
$H_{1}=\widehat{\Sigma}_{\underline{i},\underline{j}}$,
$H_{2}=\widehat{\Sigma}_{\underline{j},\underline{l}}$,
$a=x_{\underline{i},\underline{j}}$,
$b=x_{\underline{j},\underline{l}}$ and applying
Lemma~\ref{L:specialmackey}, we obtain
\[
\begin{array}{rcl}
\xi_{\underline{i},\underline{j}}\xi_{\underline{j},\underline{l}}
&=&T_{\widehat{\Sigma}_{\underline{i},\underline{j}},\widehat{\Sigma}_{r}}(x_{\underline{i},\underline{j}})T_{\widehat{\Sigma}_{\underline{j},\underline{l}},\widehat{\Sigma}_{r}}(x_{\underline{j},\underline{l}})\\
&=&\sum_{\delta\in
\widehat{\Sigma}_{\underline{j},\underline{l}}\backslash\widehat{\Sigma}_{r}/\widehat{\Sigma}_{\underline{i},\underline{j}}}T_{\widehat{\Sigma}_{\underline{i},\underline{j}}\cap\widehat{\Sigma}_{\underline{j},\underline{l}}^{\delta},\widehat{\Sigma}_{r}}(x_{\underline{i},\underline{j}}x_{\underline{j},\underline{l}}^{\delta})\\
&=&\sum_{\delta\in
\widehat{\Sigma}_{\underline{j},\underline{l}}\backslash\widehat{\Sigma}_{r}/\widehat{\Sigma}_{\underline{i},\underline{j}}}T_{\widehat{\Sigma}_{\underline{i},\underline{j}}\cap\widehat{\Sigma}_{\underline{j}\delta,\underline{l}\delta},\widehat{\Sigma}_{r}}(x_{\underline{i},\underline{j}}x_{\underline{j}\delta,\underline{l}\delta})\\
&=&\sum_{\delta\in
\widehat{\Sigma}_{\underline{j},\underline{l}}\backslash\widehat{\Sigma}_{\underline{j}}/\widehat{\Sigma}_{\underline{i},\underline{j}}}T_{\widehat{\Sigma}_{\underline{i},\underline{j},\underline{l}\delta},\widehat{\Sigma}_{r}}(x_{\underline{i},\underline{l}\delta})\\
&=&\sum_{\delta\in
\widehat{\Sigma}_{\underline{j},\underline{l}}\backslash\widehat{\Sigma}_{\underline{j}}/\widehat{\Sigma}_{\underline{i},\underline{j}}}T_{\widehat{\Sigma}_{\underline{i},\underline{l}\delta},\widehat{\Sigma}_{r}}T_{\widehat{\Sigma}_{\underline{i},\underline{j},\underline{l}\delta},\widehat{\Sigma}_{\underline{i},\underline{l}\delta}}(x_{\underline{i},\underline{l}\delta})\\
&=&\sum_{\delta\in
\widehat{\Sigma}_{\underline{j},\underline{l}}\backslash\widehat{\Sigma}_{\underline{j}}/\widehat{\Sigma}_{\underline{i},\underline{j}}}[\widehat{\Sigma}_{\underline{i},\underline{l}\delta}:\widehat{\Sigma}_{\underline{i},\underline{j},\underline{l}
\delta}]\xi_{\underline{i},\underline{l}\delta}\\
\end{array}
\]
where  the second to last equality follows by
Lemma~\ref{L:transitivity}.

The following version only involves the symmetric group
$\Sigma_{r}$, and hence easier to calculate.

\begin{Cor} For $\underline{i},\underline{j},\underline{l}\in I(n,r)$,
$\varepsilon,\varepsilon'\in\mathbb{Z}^{r}$, we have
\[\xi_{\underline{i},\underline{j}+n\varepsilon}\xi_{\underline{j},\underline{l}+n\varepsilon'}
=\sum_{\delta\in
{\Sigma}_{\underline{j},\underline{l},\varepsilon'}\backslash
{\Sigma}_{\underline{j}}/{\Sigma}_{\underline{i},\underline{j},\varepsilon}}[{\Sigma}_{\underline{i},\underline{l}\delta,\varepsilon'\delta+\varepsilon}:{\Sigma}_{\underline{i},\underline{j},\underline{l}\delta,\varepsilon'\delta,\varepsilon}]\xi_{\underline{i},\underline{l}\delta+n(\varepsilon'\delta+\varepsilon)}.\]
\end{Cor}
\begin{pf}
We have
\[\begin{array}{rcl}\xi_{\underline{i},\underline{j}+n\varepsilon}\xi_{\underline{j},\underline{l}+n\varepsilon'}
&=&
\xi_{\underline{i}-n\varepsilon,\underline{j}}\xi_{\underline{j},\underline{l}+n\varepsilon'}\\
&=&\sum_{\delta\in
\widehat{\Sigma}_{\underline{j},\underline{l}+n\varepsilon'}\backslash\widehat{\Sigma}_{\underline{j}}/\widehat{\Sigma}_{\underline{i}-n\varepsilon,\underline{j}}}[\widehat{\Sigma}_{\underline{i}-n\varepsilon,(\underline{l}+n\varepsilon')\delta}:\widehat{\Sigma}_{\underline{i}-n\varepsilon,\underline{j},(\underline{l}+n\varepsilon')
\delta}]\xi_{\underline{i}-n\varepsilon,(\underline{l}+n\varepsilon')\delta}\\
&=&\sum_{\delta\in
\Sigma_{\underline{j},\underline{l},\varepsilon'}\backslash\Sigma_{\underline{j}}/\Sigma_{\underline{i},\underline{j},\varepsilon}}[\widehat{\Sigma}_{\underline{i}-n\varepsilon,\underline{l}\delta+n\varepsilon'\delta}^{\varepsilon}:\widehat{\Sigma}_{\underline{i}-n\varepsilon,\underline{j},\underline{l}\delta+n\varepsilon'\delta}^{\varepsilon}]\xi_{\underline{i},\underline{l}\delta+n(\varepsilon'\delta+\varepsilon)}\\
&=&\sum_{\delta\in
\Sigma_{\underline{j},\underline{l},\varepsilon'}\backslash\Sigma_{\underline{j}}/\Sigma_{\underline{i},\underline{j},\varepsilon}}[\widehat{\Sigma}_{\underline{i},\underline{l}\delta+n(\varepsilon'\delta+\varepsilon)}:\widehat{\Sigma}_{\underline{i},\underline{j}+n\varepsilon,\underline{l}\delta+n(\varepsilon'\delta+\varepsilon)}]\xi_{\underline{i},\underline{l}\delta+n(\varepsilon'\delta+\varepsilon)}\\
&=&\sum_{\delta\in
\Sigma_{\underline{j},\underline{l},\varepsilon'}\backslash\Sigma_{\underline{j}}/\Sigma_{\underline{i},\underline{j},\varepsilon}}[\Sigma_{\underline{i},\underline{l}\delta,\varepsilon'\delta+\varepsilon}:\Sigma_{\underline{i},\underline{j},\varepsilon,\underline{l}\delta,\varepsilon'\delta+\varepsilon}]\xi_{\underline{i},\underline{l}\delta+n(\varepsilon'\delta+\varepsilon)}\\
&=&\sum_{\delta\in
\Sigma_{\underline{j},\underline{l},\varepsilon'}\backslash\Sigma_{\underline{j}}/\Sigma_{\underline{i},\underline{j},\varepsilon}}[\Sigma_{\underline{i},\underline{l}\delta,\varepsilon'\delta+\varepsilon}:\Sigma_{\underline{i},\underline{j},\underline{l}\delta,\varepsilon'\delta,\varepsilon}]\xi_{\underline{i},\underline{l}\delta+n(\varepsilon'\delta+\varepsilon)}\\
\end{array}\]
as desired. \hfill$\square$
\end{pf}

\vskip20pt
\section{Maps between affine Schur algebras}\label{S:relations}

In this section we will study relations between Schur algebras and
affine Schur algebras.

Let $G$ be a semigroup with identity $1_{G}$. Recall that $K^{G}$ is
an algebra and $\Delta_{G}:K^{G}\rightarrow K^{G\times G}$ and
$\epsilon_{G}:K^{G}\rightarrow K$ are algebra homomorphisms. Let $A$
be a subspace of $K^{G}$ with a fixed basis $\{a_{i}\}_{i\in I}$
such that for a fixed $g\in G$ almost all $a_{i}(g)=0$ ($i\in I$).
Assume that $\overline{A}$ is a subspace of $K^{G}$ via $(\sum_{i\in
I}\lambda_{i}a_{i})(g)=\sum_{i\in I}\lambda_{i}a_{i}(g)$. Then
$\overline{A\otimes A}$ is considered as a subspace of $ K^{G\times
G}$ via $(\sum_{i,j\in I}\lambda_{ij}a_{i}\otimes
a_{j})(g,g')=\sum_{i,j\in I}\lambda_{ij}a_{i}(g)a_{j}(g')$. For a
proof, note that for $g\in G$ the function $\sum_{i,j\in
I}\lambda_{ij}a_{i}(g)a_{j}$ lies in $K^{G}$. If $\sum_{i,j\in
I}\lambda_{ij}a_{i}\otimes a_{j}$ considered as a function on
$G\times G$ is $0$, then the coefficient $\sum_{i\in
I}\lambda_{ij}a_{i}(g)$ of $a_{j}$ equals to $0$, and hence
$\lambda_{ij}=0$.

We call $A$ a sub formal coalgebra of $K^{G}$ if in addition $A$ is
a formal coalgebra with respect to $\Delta_{G}|_{A}$ and
$\epsilon_{G}|_{A}$. Let $A$ be a sub formal coalgebra of $K^{G}$.
Define the evaluation map $e_{G}:KG\rightarrow
A^{\#},~g\mapsto(e_{G}(g):a\mapsto a(g),a\in A)$. This is a
surjective algebra homomorphism. Suppose on the contrary that
$\text{Im}(e_{G})$ is a proper subspace of $A^{\#}$. Then there
exists a nonzero element $c$ in $\overline{A}=(A^{\#})^{*}$ such
that $c(\text{Im}(e_{G}))=0$. In particular, $c(e_{G}(g))=0$ for any
$g\in G$, i.e. $c(g)=0$ for any $g\in G$. Thus $c=0$ in $K^{G}$,
contradicting the assumption that $\overline{A}$ is a subspace of
$K^{G}$.


Let $\phi:G\rightarrow H$ be a homomorphism of semigroups. Denote by
$\phi^{*}$ the algebra homomorphism from $K^{H}$ to $K^{G}$ sending
$f$ to $f\circ\phi$. If $\phi$ is injective then $\phi^{*}$ is
surjective; if $\phi$ is surjective then $\phi^{*}$ is injective.
Moreover, if $\psi:H\rightarrow L$ is also a semigroup homomorphism,
then $(\psi\circ\phi)^{*}=\phi^{*}\circ\psi^{*}$.

For $f\in K^{H}$ and $g,g'\in G$, we have
\[\begin{array}{rcl}
(\Delta_{G}\circ\phi^{*}(f))(g,g')\hspace{-7pt} & = \hspace{-7pt}&
\phi^{*}(f)(gg')=f(\phi(gg')) =
f(\phi(g)\phi(g'))\\
\hspace{-7pt}& =\hspace{-7pt} &
\Delta_{H}(f)(\phi(g),\phi(g'))=((\phi\times\phi)^{*}\circ\Delta_{H}(f))(g,g')\\
\epsilon_{G}\circ\phi^{*}(f) \hspace{-7pt}& = \hspace{-7pt}&
\phi^{*}(f)(1_{G})=f(\phi(1_{G}))=f(1_{H})=\epsilon_{H}(f).\end{array}\]
Hence $\Delta_{G}\circ
\phi^{*}=(\phi\times\phi)^{*}\circ\Delta_{H}$, and
$\epsilon_{G}\circ\phi^{*}=\epsilon_{H}$. That is, the following
diagrams commute:
\[{\setlength{\unitlength}{0.7pt}
\begin{picture}(160,90)
\put(17,60){$K^{H}$} \put(115,60){ $K^{G}$} \put(0,0){$K^{H\times
H}$} \put(100,0){ $K^{G\times G}$}
\put(28,30){\footnotesize$\Delta_{H}$}
\put(70,67){\footnotesize$\phi^{*}$}
\put(128,30){\footnotesize$\Delta_{G}$}
\put(52,9){\footnotesize$(\phi\times\phi)^{*}$}
\drawline(40,65)(105,65) \drawline(105,65)(101,61)
\drawline(105,65)(101,69) \drawline(50,5)(95,5)
\drawline(95,5)(91,1) \drawline(95,5)(91,9) \drawline(23,50)(23,20)
\drawline(23,20)(19,24) \drawline(23,20)(27,24)
\drawline(125,50)(125,20) \drawline(125,20)(121,24)
\drawline(125,20)(129,24)
\end{picture}\ \ \ \ \ \ \ \  \begin{picture}(130,90)
\put(17,60){$K^{H}$} \put(17,0){$K$} \put(75,60){ $K^{G}$}
\put(75,0){ $K$} \put(48,7){\footnotesize$\text{id}$}
\put(48,67){\footnotesize$\phi^{*}$} \drawline(35,65)(70,65)
\drawline(70,65)(66,61) \drawline(70,65)(66,69)
\drawline(22.5,50)(22.5,20) \drawline(22.5,20)(18.5,24)
\drawline(22.5,20)(26.5,24) \drawline(35,5)(70,5)
\drawline(70,5)(66,1) \drawline(70,5)(66,9) \drawline(85,50)(85,20)
\drawline(85,20)(81,24) \drawline(85,20)(89,24)
\put(28,30){$\epsilon_{H}$} \put(92,30){$\epsilon_{G}$}
\end{picture}}\]

Let $A$ and $B$ be a sub formal coalgebras of $K^{G}$ and $K^{H}$
respectively. Suppose $\phi^{*}(B)\subseteq \overline{A}$ and
$\phi^{*}|_{B}:B\rightarrow \overline{A}$ is row finite. Then
$\phi^{*}|_{B}$ is a formal homomorphism of formal coalgebras (see
Definition~\ref{D:formalhomomorphism}) and
$\overline{\phi^{*}|_{B}}=\phi^{*}|_{\overline{B}}$. By
Theorem~\ref{T:fromcoalgebratoalgebra}
$(\phi^{*}|_{B})^{\#}:A^{\#}\rightarrow B^{\#}$ is a homomorphism of
$K$-algebras. Moreover, for $g\in G$ and $b\in B$,
\[\begin{array}{rcl}
(\phi^{*}|_{B})^{\#}\circ e_{G}(g)(b)\hspace{-7pt}&=\hspace{-7pt}&
e_{G}(g)(\phi^{*}(b))=\phi^{*}(b)(g)=b(\phi(g))\\
\hspace{-7pt}&=\hspace{-7pt}&e_{H}(\phi(g))(b)=e_{H}\circ
\phi(g)(b).
\end{array}\]
Hence $(\phi^{*}|_{A})^{\#}\circ e_{G}=e_{H}\circ \phi$, i.e. the
following diagram commutes:
\[{\setlength{\unitlength}{0.7pt}
\begin{picture}(160,90)
\put(10,60){$KG$} \put(115,60){ $KH$} \put(10,0){$A^{\#}$}
\put(115,0){ $B^{\#}$} \put(28,30){\footnotesize$e_{G}$}
\put(70,67){\footnotesize$\phi$} \put(128,30){\footnotesize$e_{H}$}
\put(52,9){\footnotesize$(\phi^{*}|_{A})^{\#}$}
\drawline(40,65)(105,65) \drawline(105,65)(101,61)
\drawline(105,65)(101,69) \drawline(40,5)(105,5)
\drawline(105,5)(101,1) \drawline(105,5)(101,9)
\drawline(23,50)(23,20) \drawline(23,20)(19,24)
\drawline(23,20)(27,24) \drawline(125,50)(125,20)
\drawline(125,20)(121,24) \drawline(125,20)(129,24)
\end{picture}}\]


\subsection{The action of the extended affine Weyl group on the
affine Schur algebra}

Recall that the extended affine Weyl group $\widehat{\Sigma}_{n}$
acts on $\widetilde{GL}_{n,a}$ by conjugation (see
Lemma~\ref{L:weylaction}). In this way we regard
$w\in\widehat{\Sigma}_{n}$ as an automorphism of
$\widetilde{GL}_{n,a}$. Then the pull-back
$w^{*}:K^{\widetilde{GL}_{n,a}}\rightarrow K^{\widetilde{GL}_{n,a}}$
is an algebra automorphism, and $(ww')^{*}=w'^{*}w^{*}$.

For $i,j\in\mathbb{Z}$ and $g\in\widetilde{GL}_{n,a}$ we have
\[w^{*}(c_{ij})(g)=c_{ij}(w(g))=g'_{ij}=g_{w^{-1}(i),w^{-1}(j)}=c_{w^{-1}(i),w^{-1}(j)}(g)\]
Namely, $w^{*}(c_{ij})=c_{w^{-1}i,w^{-1}j}$. Therefore for
$\underline{i},\underline{j}\in I(\mathbb{Z},r)$ we have
\[w^{*}(c_{\underline{i},\underline{j}})=c_{w^{-1}(\underline{i}),w^{-1}(\underline{j})}.\]

Let
$f_{w}=w^{*}|_{\widetilde{A}(n,r)}:\widetilde{A}(n,r)\rightarrow\widetilde{A}(n,r)$
be the restriction of $w^{*}$ to $\widetilde{A}(n,r)$. Then $f_{w}$
is a homomorphism of formal coalgebras, $f_{w}$ is surjective, and
$\overline{f}_{w}=w^{*}|_{\overline{\widetilde{A}(n,r)}}$ is
injective. Moreover, we have $f_{w}\circ f_{w'}=f_{w'w}$. Take the
dual and we obtain an algebra automorphism
\[f_{w}^{\#}:\widetilde{S}(n,r)\rightarrow \widetilde{S}(n,r).\]
It follows by Corollary~\ref{C:formalmaps} that $f_{w'}^{\#}\circ
f_{w}^{\#}=f_{w'w}^{\#}$.

Set $w(\xi)=f_{w}^{\#}(\xi)$. This defines an action of
$\widehat{\Sigma}_{n}$ on the affine Schur algebra
$\widetilde{S}(n,r)$. Precisely, for
$\underline{i},\underline{j},\underline{p},\underline{q}\in
I(\mathbb{Z},r)$ we have
\[\begin{array}{rcl}w(\xi_{\underline{i},\underline{j}})(c_{\underline{p},\underline{q}})
\hspace{-7pt}&=\hspace{-7pt}&\xi_{\underline{i},\underline{j}}(f_{w}(c_{\underline{p},\underline{q}}))
=\xi_{\underline{i},\underline{j}}(c_{w^{-1}(\underline{p}),w^{-1}(\underline{q})})\\[5pt]
\hspace{-7pt}&=\hspace{-7pt}&\begin{cases}
1,&\text{if}~(\underline{i},\underline{j})\sim_{\widehat{\Sigma}_{r}}(w^{-1}(\underline{p}),w^{-1}(\underline{q}))\\
0,&\text{otherwise}\end{cases} \\[5pt]
\hspace{-7pt}&=\hspace{-7pt}&\begin{cases}
1,&\text{if}~(w(\underline{i}),w(\underline{j}))\sim_{\widehat{\Sigma}_{r}}(\underline{p},\underline{q})\\
0,&\text{otherwise}.\end{cases} \end{array}\] Therefore
$w(\xi_{\underline{i},\underline{j}})=\xi_{w(\underline{i}),w(\underline{j})}$.
Thus we have

\begin{Prop}\label{P:weylaction}
The $\widehat{\Sigma}_{n}$-action on $\widetilde{GL}_{n,a}$ induces
an $\widehat{\Sigma}_{n}$-action on $\widetilde{S}(n,r)$:
$w(\xi_{\underline{i},\underline{j}})=\xi_{w(\underline{i}),w(\underline{j})}$
where $w\in\widehat{\Sigma}_{n}$, $\underline{i},\underline{j}\in
I(\mathbb{Z},r)$. In particular,
$\rho(\xi_{\underline{i},\underline{j}})=\xi_{\underline{i}-(1\ldots
1),\underline{j}-(1\ldots 1)}$.\hfill$\square$
\end{Prop}

\subsection{Endomorphisms of the affine Schur algebra}

Let $a\in K^{\times},s\in\mathbb{Z}\backslash\{0\}$. Recall that
$\eta_{a,s}:\widetilde{GL}_{n,a}\rightarrow\widetilde{GL}_{n,1}$ is
an injective semigroup homomorphism (see Lemma~\ref{L:etarestricted}
(ii)). The pull-back
$\eta_{a,s}^{*}:K^{\widetilde{GL}_{n,1}}\rightarrow
K^{\widetilde{GL}_{n,a}},~~f\mapsto f\circ\eta_{a,s}$ is  an algebra
homomorphism. For $i,j=1,\ldots,n$, $l\in\mathbb{Z}$ and
$g\in\widetilde{GL}_{n,a}$, we have
\[\eta_{a,s}^{*}(c_{i,j+ln})(g)=c_{i,j+ln}(\eta_{a,s}(g))=\begin{cases} a^{\frac{l}{s}}c_{i,j+\frac{l}{s}n}(g), & \text{if}~s~|~l,\\
                                                                        0, &\text{otherwise.}
                                                                        \end{cases}\]
i.e.
\[\eta_{a,s}^{*}(c_{i,j+ln})=\begin{cases} a^{\frac{l}{s}}c_{i,j+\frac{l}{s}n}, & \text{if}~s~|~l,\\
                                                                        0, &
                                                                        \text{otherwise}.
                                                                        \end{cases}\]
For
$\varepsilon=(\varepsilon_{1},\ldots,\varepsilon_{r})\in\mathbb{Z}^{r}$
we define
$\text{ht}(\varepsilon)=\varepsilon_{1}+\ldots+\varepsilon_{r}$.
Then for $\underline{i},\underline{j}\in I(n,r)$ and
$\varepsilon\in\mathbb{Z}^{r}$ we have
\[\eta_{a,s}^{*}(c_{\underline{i},\underline{j}+n\varepsilon})=\begin{cases} a^{\text{ht}(\frac{\varepsilon}{s})}c_{\underline{i},\underline{j}+n\frac{\varepsilon}{s}}, & \text{if}~s~|~\varepsilon_{1},\ldots,\varepsilon_{r},\\
                                                                        0, &
                                                                        \text{otherwise}.
                                                                        \end{cases}\]

Let
$\varphi_{a,s}=\eta_{a,s}^{*}|_{\widetilde{A}(n,r)}:\widetilde{A}(n,r)\rightarrow\widetilde{A}(n,r)$
be the restriction of $\eta_{a,s}^{*}$ to $\widetilde{A}(n,r)$. Then
$\varphi_{a,s}$ is a surjective homomorphism of formal coalgebras.
Taking the dual we obtain an injective algebra endomorphism
\[\psi_{a,s}=\varphi_{a,s}^{\#}:\widetilde{S}(n,r)\rightarrow\widetilde{S}(n,r).\]
Precisely, for
$\underline{i},\underline{j},\underline{p},\underline{q}\in I(n,r)$
and $\varepsilon,\varepsilon'\in\mathbb{Z}^{r}$, we have
\[\begin{array}{rcl}\psi_{a,s}(\xi_{\underline{i},\underline{j}+n\varepsilon})(c_{\underline{p},\underline{q}+n\varepsilon'})
& = &
\xi_{\underline{i},\underline{j}+n\varepsilon}(\varphi_{a,s}(c_{\underline{p},\underline{q}+n\varepsilon'}))\\[5pt]
& = & \begin{cases}
\xi_{\underline{i},\underline{j}+n\varepsilon}(a^{\text{ht}(\frac{\varepsilon'}{s})}c_{\underline{p},\underline{q}+n\frac{\varepsilon'}{s}}),
& \text{if}~s~|~\varepsilon'_{1},\ldots,\varepsilon'_{r}\\
0,&\text{otherwise}\end{cases}\\[5pt]
& = & \begin{cases} a^{\text{ht}(\varepsilon)},
& \text{if}~s~|~\varepsilon'_{1},\ldots,\varepsilon'_{r},~\text{and}~(\underline{i},\underline{j}+n\varepsilon)\sim_{\Sigma_{r}}(\underline{p},\underline{q}+n\frac{\varepsilon'}{s})\\
0,&\text{otherwise.}\end{cases}
\end{array}\]
Therefore
$\psi_{a,s}(\xi_{\underline{i},\underline{j}+n\varepsilon})=a^{\text{ht}(\varepsilon)}\xi_{\underline{i},\underline{j}+ns\varepsilon}$.
It is easy to see that the restriction of $\psi_{a,s}$ to $S(n,r)$
is the identity map. \vskip10pt

Recall that $\eta_{a}:\widetilde{GL}_{n,a}\rightarrow GL_{n}$ is a
surjective semigroup homomorphism (see Lemma~\ref{L:etarestricted}
(ii)). The pull-back $\eta_{a}^{*}:K^{GL_{n}}\rightarrow
K^{\widetilde{GL}_{n,a}},~~f\mapsto f\circ\eta_{a}$ is an injective
algebra homomorphism. For $i,j=1,\ldots,n$ and
$g\in\widetilde{GL}_{n,a}$ we have
\[\eta_{a}^{*}(c_{ij})(g)=c_{ij}(\eta_{a}(g))=\sum_{l\in\mathbb{Z}}a^{l}g_{i,j+ln}=\sum_{l\in\mathbb{Z}}a^{l}c_{i,j+ln}(g)\]
i.e. $\eta_{a}^{*}(c_{ij})=\sum_{l\in\mathbb{Z}}a^{l}c_{i,j+ln}$.
Therefore for $\underline{i},\underline{j}\in I(n,r)$ we have
\[\eta_{a}^{*}(c_{\underline{i},\underline{j}})=\sum_{\varepsilon\in\mathbb{Z}^{r}}a^{\text{ht}(\varepsilon)}c_{\underline{i},\underline{j}+n\varepsilon}.\]

Let $\varphi_{a}:A(n,r)\rightarrow\overline{\widetilde{A}(n,r)}$
denote the restriction of $\eta_{a}^{*}$ to $A(n,r)$. Then
$\varphi_{a}$ is a formal homomorphism of formal coalgebras.
Moreover $\overline{\varphi_{a}}=\varphi_{a}$ is injective. Taking
the dual we obtain a surjective algebra homomorphism
\[\psi_{a}=\varphi_{a}^{\#}:\widetilde{S}(n,r)\rightarrow S(n,r).\]
Precisely, for
$\underline{i},\underline{j},\underline{p},\underline{q}\in I(n,r)$
and $\varepsilon\in\mathbb{Z}^{r}$, we have
\[\begin{array}{rcl}\psi_{a}(\xi_{\underline{i},\underline{j}+n\varepsilon})(c_{\underline{p},\underline{q}})
& = &
\xi_{\underline{i},\underline{j}+n\varepsilon}(\varphi_{a}(c_{\underline{p},\underline{q}}))\\
& = &
\xi_{\underline{i},\underline{j}+n\varepsilon}(\sum_{\varepsilon'\in\mathbb{Z}^{r}}a^{\text{ht}(\varepsilon')}c_{\underline{p},\underline{q}+n\varepsilon'})\\
& = &
\sum_{\varepsilon'\in\mathbb{Z}^{r}}a^{\text{ht}(\varepsilon')}\xi_{\underline{i},\underline{j}+n\varepsilon}(c_{\underline{p},\underline{q}+n\varepsilon'})\\
& = & \begin{cases}
a^{\text{ht}(\varepsilon)}[\Sigma_{\underline{i},\underline{j}}:\Sigma_{\underline{i},\underline{j},\varepsilon}],
& \text{if}~(\underline{i},\underline{j})\sim_{\Sigma_{r}}(\underline{p},\underline{q})\\
0,&\text{otherwise}.\end{cases}
\end{array}\]
Therefore~$\psi_{a}(\xi_{\underline{i},\underline{j}+n\varepsilon})=a^{\text{ht}(\varepsilon)}[\Sigma_{\underline{i},\underline{j}}:\Sigma_{\underline{i},\underline{j},\varepsilon}]\xi_{\underline{i},\underline{j}}$.
It is easy to see that the restriction of $\psi_{a}$ to $S(n,r)$ is
the identity map.

In summary we have
\begin{Prop}\label{P:psi} Let $a\in K^{\times}$, $s\in\mathbb{Z}\backslash\{0\}$.

{\rm (i)} The injective semigroup homomorphism $\eta_{a,s}$ induces
an injective algebra endomorphism of the affine Schur algebra
\[\psi_{a,s}:\widetilde{S}(n,r)\rightarrow \widetilde{S}(n,r),~~\xi_{\underline{i},\underline{j}+n\varepsilon}\mapsto a^{\rm{ht}(\varepsilon)}\xi_{\underline{i},\underline{j}+ns\varepsilon}\]
where $\underline{i},\underline{j}\in I(n,r)$ and
$\varepsilon\in\mathbb{Z}^{r}$.

{\rm (ii)} The surjective semigroup homomorphism $\eta_{a}$ induces
a surjective algebra homomorphism
\[\psi_{a}:\widetilde{S}(n,r)\rightarrow S(n,r),~~\xi_{\underline{i},\underline{j}+n\varepsilon}\mapsto a^{\rm{ht}(\varepsilon)}[\Sigma_{\underline{i},\underline{j}}:\Sigma_{\underline{i},\underline{j},\varepsilon}]\xi_{\underline{i},\underline{j}}\]
where $\underline{i},\underline{j}\in I(n,r)$ and
$\varepsilon\in\mathbb{Z}^{r}$.

{\rm (iii)} The restrictions of $\psi_{a,s}$ and $\psi_{a}$ to
$S(n,r)$ are the identity map.\hfill$\square$
\end{Prop}

Denote by $\psi_{a,0}$ the composition
$\widetilde{S}(n,r)\stackrel{\psi_{a}}{\rightarrow}S(n,r)
\hookrightarrow\widetilde{S}(n,r)$. Then for any $s\in\mathbb{Z}$ we
have $\psi_{a,s}(\xi_{\underline{i},\underline{j}+n\varepsilon})=
a^{\text{ht}(\varepsilon)}[\Sigma_{\underline{i},\underline{j}}:\Sigma_{\underline{i},\underline{j},\varepsilon}]^{\delta_{s,0}}\xi_{\underline{i},\underline{j}+ns\varepsilon}$.

\begin{Prop}\label{P:psi} Assume $a,a'\in K^{\times}$ and
$s,s'\in\mathbb{Z}$. Then
$\psi_{a,s}\circ\psi_{a',s'}=\psi_{a'a^{s'},ss'}$.
\end{Prop}

\begin{pf} For $\underline{i},\underline{j}\in
I(n,r)$ and $\varepsilon\in\mathbb{Z}^{r}$ we have
\[\begin{array}{rcl}
\psi_{a,s}\circ\psi_{a',s'}(\xi_{\underline{i},\underline{j}+n\varepsilon})
&=&\psi_{a,s}((a')^{\text{ht}(\varepsilon)}[\Sigma_{\underline{i},\underline{j}}:\Sigma_{\underline{i},\underline{j},\varepsilon}]^{\delta_{s',0}}\xi_{\underline{i},\underline{j}+ns'\varepsilon})\\
&=&(a')^{\text{ht}(\varepsilon)}[\Sigma_{\underline{i},\underline{j}}:\Sigma_{\underline{i},\underline{j},\varepsilon}]^{\delta_{s',0}}a^{\text{ht}(s'\varepsilon)}[\Sigma_{\underline{i},\underline{j}}:\Sigma_{\underline{i},\underline{j},s'\varepsilon}]^{\delta_{s,0}}\xi_{\underline{i},\underline{j}+nss'\varepsilon}\\
&=&(a'a^{s'})^{\text{ht}(\varepsilon)}[\Sigma_{\underline{i},\underline{j}}:\Sigma_{\underline{i},\underline{j},\varepsilon}]^{\delta_{ss',0}}\xi_{\underline{i},\underline{j}+nss'\varepsilon}\\
&=&\psi_{a'a^{s'},ss'}(\xi_{\underline{i},\underline{j}+n\varepsilon}).
\end{array}\]
Therefore $\psi_{a,s}\circ\psi_{a',s'}=\psi_{a'a^{s'},ss'}$.
\hfill$\square$
\end{pf}

\vskip10pt

The transpose $g^{tr}$ of a matrix $g\in\mathfrak{M}_{n}$ is also a
matrix in $\mathfrak{M}_{n}$. It follows by Lemma~\ref{L:eta}~(ii)
that
\[\begin{array}{c}\widetilde{\text{det}}_{a^{-1}}(g^{tr})=\text{det}(\eta_{a^{-1}}(g^{tr}))=\text{det}(\eta_{a}(g)^{tr})
=\text{det}(\eta_{a}(g))=\widetilde{\text{det}}_{a}(g).\end{array}\]
Therefore $g\in\widetilde{GL}_{n,a}$ if and only if
$g^{tr}\in\widetilde{GL}_{n,a^{-1}}$. In particular, the map
$T:g\mapsto g^{tr}$ is a semigroup anti-isomorphism from
$\widetilde{GL}_{n,a}$ to $\widetilde{GL}_{n,a^{-1}}$. Therefore
taking transpose induces an algebra isomorphism
\[T^{*}:K^{\widetilde{GL}_{n,a^{-1}}}\rightarrow K^{\widetilde{GL}_{n,a}},~~f\mapsto f\circ T.\]
For $i,j\in\mathbb{Z}$ and $g\in\widetilde{GL}_{n,a}$ we have
$T^{*}(c_{ij})(g)=c_{ij}(g^{tr})=c_{ji}(g)$, i.e.
$T^{*}(c_{ij})=c_{ji}$. Therefore for
$\underline{i},\underline{j}\in I(\mathbb{Z},r)$ we have
\[T^{*}(c_{\underline{i},\underline{j}})=c_{\underline{j},\underline{i}}.\]
The restriction map
$T^{*}|_{\widetilde{A}(n,r)}:\widetilde{A}(n,r)\rightarrow
\widetilde{A}(n,r)$ is an anti-homomorphism of formal coalgebras.
Moreover, the map
$\overline{T^{*}|_{\widetilde{A}(n,r)}}=T^{*}|_{\overline{\widetilde{A}(n,r)}}$
is injective and $T^{*}|_{\widetilde{A}(n,r)}$ is surjective. Take
the dual and we obtain an algebra anti-automorphism
\[\widetilde{J}=(T^{*}|_{\widetilde{A}(n,r)})^{\#}:\widetilde{S}(n,r)\rightarrow \widetilde{S}(n,r),~~\xi_{\underline{i},\underline{j}}\mapsto\xi_{\underline{j},\underline{i}}\]
where $\underline{i},\underline{j}\in I(\mathbb{Z},r)$. The
quantized $\widetilde{J}$ is given in~\cite{L} Lemma 1.11.

\subsection{Transfer maps}

The multiplicative function $\widetilde{\text{det}}_{a}$ on
$\widetilde{GL}_{n,a}$ can be written as
\[\widetilde{\text{det}}_{a}=\sum_{\sigma\in\Sigma_{n}}{{\text{sgn}(\sigma)}\sum_{\varepsilon\in
\mathbb{Z}^{n}}a^{\text{ht}(\varepsilon)}c_{(12\ldots n),(12\ldots
n)\sigma+n\varepsilon}}\in\overline{\widetilde{A}(n,n)}.\] For
$g,h\in~\widetilde{GL}_{n,a}$ we have
\[\begin{array}{c}\Delta(\widetilde{\text{det}}_{a})(g,h)=\widetilde{\text{det}}_{a}(gh)=\widetilde{\text{det}}_{a}(g)\widetilde{\text{det}}_{a}(h)
=(\widetilde{\text{det}}_{a}\otimes\widetilde{\text{det}}_{a})(g,h).
\end{array}\]
So
$\Delta(\widetilde{\text{det}}_{a})=\widetilde{\text{det}}_{a}\otimes
\widetilde{\text{det}}_{a}$. Therefore multiplying by
$\widetilde{\text{det}}_{a}$ is a formal homomorphism
$\widetilde{\text{det}}_{a}:\widetilde{A}(n,r)\rightarrow
\overline{\widetilde{A}(n,n+r)}$ of formal coalgebras.

\begin{Lem}\label{L:transfersimilar} We have
$\widetilde{\rm{det}}_{a}\circ\varphi_{a,1}=\overline{\varphi_{a,1}}\circ\widetilde{\rm{det}}_{1}$.
\end{Lem}
\begin{pf}
Note that
$\widetilde{\text{det}}_{a}=\overline{\varphi_{a,1}}(\widetilde{\text{det}}_{1})$.
For $c\in\widetilde{A}(n,r)$ we have
\[\begin{array}{rcl}
\widetilde{\text{det}}_{a}\circ\varphi_{a,1}(c)\hspace{-7pt}&=\hspace{-7pt}&\varphi_{a,1}(c)\widetilde{\text{det}}_{a}=\varphi_{a,1}(c)\overline{\varphi_{a,1}}(\widetilde{\text{det}}_{1})\\
\hspace{-7pt}&=\hspace{-7pt}&\overline{\varphi_{a,1}}(c~\widetilde{\text{det}}_{1})=\overline{\varphi_{a,1}}\circ\widetilde{\text{det}}_{1}(c).
\end{array}\]
Therefore
$\widetilde{\text{det}}_{a}\circ\varphi_{a,1}=\overline{\varphi_{a,1}}\circ\widetilde{\text{det}}_{1}$.\hfill$\square$
\end{pf}

It follows by Proposition~\ref{P:algindep2} that the map
$\overline{\widetilde{\text{det}}_{a}}$ is injective. Take the dual
of the map $\widetilde{\text{det}}_{a}$ and we obtain a surjective
algebra homomorphism
\[\widetilde{\text{det}}_{a}^{\#}:\widetilde{S}(n,n+r)\rightarrow \widetilde{S}(n,r),~~\xi_{\underline{i},\underline{j}}\mapsto \Big (c\mapsto
\xi_{\underline{i},\underline{j}}(c~\widetilde{\text{det}}_{a})
 {\ \rm for\ } c\in\widetilde{A}(n,r)\Big).\]

\begin{Prop}\label{T:affinedet}

{\rm (i)} We have
$\psi_{a,1}\circ\widetilde{\rm{det}}_{a}^{\#}=\widetilde{\rm{det}}_{1}^{\#}\circ\psi_{a,1}$.

{\rm (ii)} We have
$\widetilde{\rm{det}}_{a}|_{S(n,n+r)}=\rm{det}^{*}$.

{\rm (iii)} On $K\widetilde{SL}_{n,a}$ we have $\widetilde{e}_{a}^{\
r}=\widetilde{\rm{det}}_{a}^{\#}\circ \widetilde{e}_{a}^{\ n+r}$.

\end{Prop}

\begin{pf}
(i) This follows from Lemma~\ref{L:transfersimilar} and
Theorem~\ref{T:formalmaps} (v).

(ii) Let $\underline{i},\underline{j}\in I(n,n+r)$,
$\underline{p},\underline{q}\in I(n,r)$. Then
\[\begin{array}{rcl}\widetilde{\text{det}}_{a}|_{S(n,n+r)}(\xi_{\underline{i},\underline{j}})(c_{\underline{p},\underline{q}})
&=&\xi_{\underline{i},\underline{j}}(\widetilde{\text{det}}_{a}c_{\underline{p},\underline{q}})\\
&=&\xi_{\underline{i},\underline{j}}(\sum_{\sigma\in\Sigma_{n}}{{\text{sgn}(\sigma)}\sum_{\varepsilon\in
\mathbb{Z}^{n}}a^{\text{ht}(\varepsilon)}c_{(12\ldots n),(12\ldots
n)\sigma+n\varepsilon}}c_{\underline{p},\underline{q}})\\
&=&\sum_{\sigma\in\Sigma_{n}}{{\text{sgn}(\sigma)}\sum_{\varepsilon\in
\mathbb{Z}^{n}}a^{\text{ht}(\varepsilon)}\xi_{\underline{i},\underline{j}}(c_{(12\ldots
n),(12\ldots
n)\sigma+n\varepsilon}}c_{\underline{p},\underline{q}}).
\end{array}\]
Since $\underline{i},\underline{j}\in
I(n,n+r),\underline{p},\underline{q}\in I(n,r)$, we have
$\varepsilon\neq 0$ implies
$\xi_{\underline{i},\underline{j}}(c_{(12\ldots n),(12\ldots
n)\sigma+n\varepsilon}c_{\underline{p},\underline{q}})=0$. Therefore
\[\begin{array}{rcl}
\widetilde{\text{det}}_{a}|_{S(n,n+r)}(\xi_{\underline{i},\underline{j}})(c_{\underline{p},\underline{q}})
\hspace{-7pt}&=\hspace{-7pt}&\sum_{\sigma\in\Sigma_{n}}{\text{sgn}(\sigma)}\xi_{\underline{i},\underline{j}}(c_{(12\ldots
n),(12\ldots n)\sigma}c_{\underline{p},\underline{q}})\\
\hspace{-7pt}&=\hspace{-7pt}&\xi_{\underline{i},\underline{j}}(\sum_{\sigma\in\Sigma_{n}}{\text{sgn}(\sigma)}c_{(12\ldots
n),(12\ldots
n)\sigma}c_{\underline{p},\underline{q}})\\
\hspace{-7pt}&=\hspace{-7pt}&\xi_{\underline{i},\underline{j}}(\text{det}~
c_{\underline{p},\underline{q}})
=\text{det}^{*}(\xi_{\underline{i},\underline{j}})(c_{\underline{p},\underline{q}}).
\end{array}\]

(iii) For $g\in \widetilde{SL}_{n,a}$, $c\in\widetilde{A}(n,r)$,
\[\begin{array}{rcl}
\widetilde{\text{det}}_{a}^{\#}\circ\widetilde{e}_{a}^{\
n+r}(g)(c)\hspace{-7pt}&=\hspace{-7pt}&\widetilde{e}_{a}^{\
n+r}(g)(c~\widetilde{\text{det}}_{a})=(c~\widetilde{\text{det}}_{a})(g)\\
\hspace{-7pt}&=\hspace{-7pt}&c(g)\widetilde{\text{det}}_{a}(g)=c(g)=\widetilde{e}_{a}^{\
r}(g)(c).
\end{array}\]
\hfill$\square$
\end{pf}

\begin{Thm}\label{T:gpandalg}
The following diagram commutes:
\[
{\setlength{\unitlength}{0.7pt}
\begin{picture}(480,200)
\put(50,130){$S(n,n+r)$}

\put(125,134){ $\line(1,0){145}$}

\drawline(277,134)(272,139) \drawline(277,134)(272,129)
\drawline(272,134)(267,139) \drawline(272,134)(267,129)

\put(280,130){ $S(n,r) $}

\put(50,0){$\widetilde{S}(n,n+r) $}

\put(125,4){ $\line(1,0){145}$}

\drawline(277,4)(272,9) \drawline(277,4)(272,-1)
\drawline(272,4)(267,9) \drawline(272,4)(267,-1)

\put(280,0){ $\widetilde{S}(n,r)$}

\put(75,115){$\line(0,-1){95}$}

\drawline(75.5,115)(70.5,110) \drawline(75.5,115)(80.5,110)
\drawline(75.5,110)(70.5,105) \drawline(75.5,110)(80.5,105)

\put(300,115){$\line(0,-1){95}$}

\drawline(300,115)(295,110) \drawline(300,115)(305,110)
\drawline(300,110)(295,105) \drawline(300,110)(305,105)

\put (170,30){$K\widetilde{SL}_{n,a}$}

\put(170,160){$KSL_{n}$}

\drawline(100,146)(160,159)

 \drawline(100,146)(104,151) \drawline(100,146)(105,143)
 \drawline(105,147)(109,152) \drawline(105,147)(110,144)

\drawline(230,160)(280,146)

\drawline(280,146)(277,152) \drawline(280,146)(275,143)
\drawline(275,147)(272,153) \drawline(275,147)(270,144)

\dashline{5}(100,16)(165,29)

 \drawline(100,16)(104,21) \drawline(100,16)(105,13)
 \drawline(105,17)(109,22) \drawline(105,17)(110,14)

\dashline{5}(230,30)(280,16)

\drawline(280,16)(277,22) \drawline(280,16)(275,13)
\drawline(275,17)(272,23) \drawline(275,17)(270,14)

\dashline{5}(185,150)(185,45)

\drawline(185,150)(180,145) \drawline(185,150)(190,145)
\drawline(185,145)(180,140) \drawline(185,145)(190,140)

\put(370,160){$KGL_{n}$}

\drawline(365,159)(310,146)
\drawline(310,146)(314,151)\drawline(310,146)(315,143)
\drawline(315,147)(319,152) \drawline(315,147)(320,144)

\drawline(230,168)(355,168) \put(235,170.5){$\oval(20,5)[bl]$}
\put(235,174){$\oval(20,5)[tl]$} \drawline(355,168)(350,173)
\drawline(355,168)(350,163)

\put(370,30){$K\widetilde{GL}_{n,a}$}

\drawline(365,29)(310,16)
\drawline(310,16)(314,21)\drawline(310,16)(315,13)
\drawline(315,17)(319,22) \drawline(315,17)(320,14)

\dashline{5}(235,38)(353,38) \put(240,40.5){$\oval(20,5)[bl]$}
\put(240,44){$\oval(20,5)[tl]$} \drawline(355,38)(350,43)
\drawline(355,38)(350,33)

\put(385,145){$\line(0,-1){95}$}
 \drawline(385.5,145)(380.5,140) \drawline(385.5,145)(389.5,140)
 \drawline(385.5,140)(380.5,135) \drawline(385.5,140)(389.5,135)

\put(80,60){$\psi_{a}$} \put(195,90){$\eta_{a}$}
\put(305,60){$\psi_{a}$} \put(395,90){$\eta_{a}$}

\put(130,160){$e^{n+r}$} \put(240,143){$e^{r}$}
\put(130,30){$\widetilde{e}_{a}^{\ n+r}$}
\put(233,13){$\widetilde{e}_{a}^{\ r}$}
\put(150,120){$\rm{det}^{*}$}
\put(150,-15){$\widetilde{\rm{det}}_{a}^{\#}$}
\end{picture}}
\]
\end{Thm}
\vskip20pt
\begin{pf}
By Proposition~\ref{T:affinedet} it suffices  to prove
$\psi_{a}\circ\widetilde{\text{det}}_{a}^{\#}=\text{det}^{*}\circ\psi_{a}$.
Following the commutative parts of the diagram we have
\[\begin{array}{rcl}\psi_{a}\circ\widetilde{\text{det}}_{a}^{\#}\circ\widetilde{e}_{a}^{~n+r}
\hspace{-7pt}&=\hspace{-7pt}&\psi_{a}\circ\widetilde{e}_{a}^{~r}=e^{r}\circ\eta_{a}\\
\hspace{-7pt}&=\hspace{-7pt}&\text{det}^{*}\circ e^{n+r}\circ
\eta_{a}=\text{det}^{*}\circ\psi_{a}\circ \widetilde{e}_{a}^{~n+r}.
\end{array}\]
It follows from the surjectivity of $\widetilde{e}_{a}^{~n+r}$ that
$\psi_{a}\circ\widetilde{\text{det}}_{a}^{\#}=\text{det}^{*}\circ\psi_{a}$.
\hfill$\square$
\end{pf}

\section{The Lie algebras $\mathfrak{g}_{n}$ and
$\mathfrak{g}'_{n}$}\label{S:liealgebras}

In this section, $K=\mathbb{C}$, and $n\geq 2$.

Let $\mathfrak{g}_{n}$ be the underlying Lie algebra of
$\mathfrak{M}_{n}$ with Lie bracket the commutator, and
$\mathfrak{g}'_{n}=[\mathfrak{g}_{n},\mathfrak{g}_{n}]$ its Lie
subalgebra. Then $\mathfrak{g}_{n}=\mathfrak{gl}_{n}[t,t^{-1}]$ and
$\mathfrak{g}'_{n}=\mathfrak{sl}_{n}[t,t^{-1}]$ are loop algebras,
and they are quotients of the affine general linear Lie algebra
$\widehat{\mathfrak{gl}}_{n}$ and the affine special linear Lie
algebra $\widehat{\mathfrak{sl}}_{n}$ respectively.

\begin{Lem}(\cite{Ka}) As a Lie algebra over $\mathbb{C}$, the loop algebra $\mathfrak{g}'_{n}$ is generated by
$\{E_{s,s+1},E_{s,s-1}|s=1,\ldots,n\}$.\hfill$\square$
\end{Lem}

The universal enveloping algebra $\mathcal{U}(\mathfrak{g}_{n})$ of
$\mathfrak{g}_{n}$ acts naturally on $\widetilde{E}$, and hence on
$\widetilde{E}^{\otimes r}$ via the comultiplication. Let
$\widetilde{\pi}:\mathcal{U}(\mathfrak{g}_{n})\rightarrow
\text{End}_{\mathbb{C}}(\widetilde{E}^{\otimes r})$ be the
corresponding representation map.

\begin{Lem}
This action commutes with the right action of
$\widehat{\Sigma}_{r}$. In particular, the image ${\rm
Im}(\widetilde{\pi})$ of $\widetilde{\pi}$ is a subalgebra of
$\widetilde{S}(n,r)$.
\end{Lem}
\begin{pf}

Let $s\in\{1,\ldots,n\}$, $t\in\mathbb{Z}$. For $\underline{i}\in
I(\mathbb{Z},r)$ and $k=1,\ldots,r$ we define
$\underline{i}^{k}=(i_{1},\ldots,i_{k-1},s-t+i_{k},i_{k+1},\ldots,i_{r})$.

Let $\underline{i}\in I(\mathbb{Z},r)$, and
$w=(\sigma,\varepsilon)\in\widehat{\Sigma}_{r}$ with
$\sigma\in\Sigma_{r}$, $\varepsilon\in\mathbb{Z}^{r}$. We have
\[
\begin{array}{rcl}
(E_{st}v_{\underline{i}})w
\hspace{-7pt}&=\hspace{-7pt}&((\sum_{k=1}^{r}1^{\otimes k-1}\otimes
E_{st}\otimes 1^{\otimes r-k})v_{\underline{i}})w\\
\hspace{-7pt}&=\hspace{-7pt}&(\sum_{k=1}^{r}\delta_{\overline{t},\overline{i}_{k}}v_{\underline{i}^{k}})w
=\sum_{k=1}^{r}\delta_{\overline{t},\overline{i}_{k}}v_{\underline{i}^{k}w}\\[8pt]
E_{st}(v_{\underline{i}}w)
\hspace{-7pt}&=\hspace{-7pt}&E_{st}(v_{\underline{i}w})
=\sum_{k=1}^{r}\delta_{\overline{t},\overline{(\underline{i}w)}_{k}}v_{(\underline{i}w)^{k}}\\
\hspace{-7pt}&=\hspace{-7pt}&\sum_{k=1}^{r}\delta_{\overline{t},\overline{(\underline{i}\sigma)}_{k}}v_{(\underline{i}\sigma+n\varepsilon)^{k}}
=\sum_{k=1}^{r}\delta_{\overline{t},\overline{i}_{\sigma^{-1}(k)}}v_{(\underline{i}\sigma)^{k}+n\varepsilon}\\
\hspace{-7pt}&=\hspace{-7pt}&\sum_{k=1}^{r}\delta_{\overline{t},\overline{i}_{\sigma^{-1}(k)}}v_{\underline{i}^{\sigma^{-1}(k)}\sigma+n\varepsilon}\\
\hspace{-7pt}&=\hspace{-7pt}&\sum_{k=1}^{r}\delta_{\overline{t},\overline{i}_{k}}v_{\underline{i}^{k}\sigma+n\varepsilon}
=\sum_{k=1}^{r}\delta_{\overline{t},\overline{i}_{k}}v_{\underline{i}^{k}w}.
\end{array}
\]
This completes the proof. \hfill$\square$
\end{pf}

\begin{Lem}\label{L:imageofgen} Let $s=1,\ldots,n$, $t\in\mathbb{Z}$, $s\not= t$. Then
\[\begin{array}{rcl}
\widetilde{\pi}(E_{st})&=&\sum_{\underline{i}\in
I(n,r-1)/\Sigma_{r-1}}\xi_{\underline{i}s, \underline{i}t}\\
\widetilde{\pi}(E_{ss})&=&\sum_{\underline{i}\in
I(n,r)/\Sigma_{r}}\lambda_{\underline{i},s}\xi_{\underline{i},\underline{i}}\end{array}\]
where $\lambda_{\underline{i},s}$ is the number of $s$ in
$\underline{i}$.
\end{Lem}
\begin{pf} Let $\underline{q}\in
I(\mathbb{Z},r)$. We have
\[
\begin{array}{rcl}
(\sum_{\underline{i}}\xi_{\underline{i}s,\underline{i}t})(v_{\underline{q}})
\hspace{-7pt}&=\hspace{-7pt}&\sum_{\underline{i}}\xi_{\underline{i}s,\underline{i}t}(v_{\underline{q}})
=\sum_{\underline{i}}\sum_{\underline{p}\in
I(\mathbb{Z},r)}\xi_{\underline{i}s,\underline{i}t}(c_{\underline{p},\underline{q}})v_{\underline{p}}\\
\hspace{-7pt}&=\hspace{-7pt}&\sum_{\underline{i}}\sum_{k=1}^{r}\xi_{\underline{i}s,\underline{i}t}(c_{\underline{q}^{k},\underline{q}})v_{\underline{q}^{k}}
=\sum_{k=1}^{r}\delta_{\overline{t},\overline{q}_{k}}v_{\underline{q}^{k}}\\
\hspace{-7pt}&=\hspace{-7pt}&E_{st}(v_{\underline{q}})\\[8pt]
(\sum_{\underline{i}}\lambda_{\underline{i},s}\xi_{\underline{i},\underline{i}})(v_{\underline{q}})
\hspace{-7pt}&=\hspace{-7pt}&\sum_{\underline{i}}\lambda_{\underline{i},s}\xi_{\underline{i},\underline{i}}(v_{\underline{q}})
=\lambda_{\underline{q}}v_{\underline{q}}
=\sum_{k=1}^{r}\delta_{s,\overline{q}_{k}}v_{\underline{q}^{k}}\\
\hspace{-7pt}&=\hspace{-7pt}&E_{ss}(v_{\underline{q}})
\end{array}
\]
as desired. \hfill$\square$
\end{pf}

\begin{Lem}\label{L:imageofgen-gln}
The restriction map
$\widetilde{\pi}|_{\mathcal{U}(\mathfrak{gl}_{n})}$ equals to $\pi$
(defined in Section~\ref{S:gln}).
\end{Lem}
\begin{pf} It is enough to prove that if $s,t=1,\ldots,n$, $s\not= t$, then
\[\begin{array}{rcl}
\pi(E_{st})&=&\sum_{\underline{i}\in
I(n,r-1)/\Sigma_{r-1}}\xi_{\underline{i}s, \underline{i}t}\\
\pi(E_{ss})&=&\sum_{\underline{i}\in
I(n,r)/\Sigma_{r}}\lambda_{\underline{i},s}\xi_{\underline{i},\underline{i}}.\end{array}\]
The proof is the same as that of Lemma~\ref{L:imageofgen} except
that we need to replace $\underline{q}\in I(\mathbb{Z},r)$ by
$\underline{q}\in I(n,r)$. \hfill$\square$
\end{pf}

By Lemma~\ref{L:imageofgen-gln}
$\widetilde{\pi}(\mathcal{U}(\mathfrak{sl}_{n}))=\pi(\mathcal{U}(\mathfrak{sl}_{n}))=S(n,r)$
contains all $\xi_{\underline{i},\underline{i}}$, $\underline{i}\in
I(n,r)$. Therefore by Lemma~\ref{L:imageofgen}
$\widetilde{\pi}(\mathcal{U}(\mathfrak{g}_{n}))$ contains a
subalgebra of $\widetilde{S}(n,r)$ generated by $Y$
 and $\widetilde{\pi}(\mathcal{U}(\mathfrak{g}'_{n}))$ contains a subalgebra of $\widetilde{S}(n,r)$  generated
 by $X=X_{1}\cup X_{2}$, where
\[\begin{array}{rcl}
Y&=&\{\xi_{\underline{i}s, \underline{i}t}|\underline{i}\in
I(n,r-1),s=1,\ldots,n,t\in\mathbb{Z}\}\\
X_{1}&=&\{\xi_{\underline{i}s, \underline{i}(s+1)}|\underline{i}\in
I(n,r-1), s=1,\ldots,n\}\\
X_{2}&=&\{\xi_{\underline{i}s, \underline{i}(s-1)}|\underline{i}\in
I(n,r-1), s=1,\ldots,n\}.\end{array}\]

\begin{Lem}\label{L:generator}
{\rm (i)} As a $K$-algebra, $\widetilde{S}(n,r)$ is generated by
$Y$.

{\rm (ii)} Assume $r<n$. As a $K$-algebra, $\widetilde{S}(n,r)$ is
generated by $X$.
\end{Lem}
\begin{pf}
(i) For $\xi_{\underline{i},\underline{j}}\in\widetilde{S}(n,r)$
with $\underline{i}\in I(n,r)$, define its index to be the number of
$s$ in $\{1,\ldots,r\}$ such that $j_{s}\not= i_{s}$. Induct on the
index. Clearly $Y$ is the set of all
$\xi_{\underline{i},\underline{j}}$'s of index $0$ and $1$.

Suppose $\xi_{\underline{i},\underline{j}}$ is of index $m\geq 2$.
Without loss of generality we may assume that $i_{1}\neq j_{1}$.
Assume $1\leq p\leq q\leq r$ are such that
\[ \begin{array}{rcl}
i_{s}\hspace{-7pt}&=\hspace{-7pt}&j_{s}=i_{1}\ {\rm for}\ 2\leq s\leq p,\\
j_{s}\hspace{-7pt}&\equiv\hspace{-7pt}& i_{1}\hspace{-7pt}\pmod{n},\ j_{s}\not= i_{s}\ {\rm for}\ p+1\leq s\leq q,\\
j_{s}\hspace{-7pt}&\not\equiv\hspace{-7pt}&
i_{1}\hspace{-7pt}\pmod{n}\ {\rm for}\ q+1\leq s\leq r.\end{array}
\]

Let $\underline{j}'=i_{1}j_{2}\ldots j_{r}$. Then
\[
\xi_{\underline{i},\underline{j}'}\xi_{\underline{j}',\underline{j}}=a_{\underline{i},\underline{j}}\xi_{\underline{i},\underline{j}}+\sum_{\underline{k}}{a_{\underline{i},\underline{k}}\xi_{\underline{i},\underline{k}}}
\]
where $\underline{k}$ \ lies in the set $\{(i_{1}\ldots
i_{1}j_{p+1}\ldots j_{s-1} (j_{s}+j_{1}-i_{1})j_{s+1}\ldots
j_{r})|p+1\leq s\leq q\}$ and $a_{\underline{i},\underline{j}}\neq
0$.

In the above equality, the elements
$\xi_{\underline{i},\underline{j}'}$,
$\xi_{\underline{j}',\underline{j}}$, and all
$\xi_{\underline{i},\underline{k}}$ have indices smaller than $m$.
So by induction we can finish the proof.

(ii) We only need to show $X$ generates $Y$. Let
$s\in\{1,\ldots,n\}$, $t\in \mathbb{Z}$. We may assume that $s\leq
t$.

Note that by Proposition~\ref{P:weylaction} the action of
$\rho\in\widehat{\Sigma}_{n}$ on $\widetilde{S}(n,r)$ permutes
$X_{1}$. So if $s\leq t<s+n$, then
$\rho^{-s+1}(\xi_{\underline{i}s,\underline{i}t})=\xi_{\underline{i}'1,\underline{i}'(t-s+1)}$
lies in $ S(n,r)$, where
$\underline{i}'=\overline{\underline{i}-(s-1)(1\ldots 1)}\in
I(n,r-1)$. Now $\xi_{\underline{i}'1,\underline{i}'(t-s+1)}$ is
generated by elements in $X_{1}\cap S(n,r)$, and hence
$\xi_{\underline{i}s,\underline{i}t}$ is generated by elements in
$\rho^{s-1}(X_{1}\cap S(n,r))\subseteq X_{1}$.

If $t\geq s+n$, then there exists $s'$ such that $s<s'<s+n$ and
\[\begin{array}{c}
s'\not\equiv i_{p} \hspace{-7pt}\pmod{n}\ {\rm for\ any}\ 1\leq
p\leq r-1.
\end{array}
\]
Therefore
\[
\xi_{\underline{i}s,\underline{i}s'}\xi_{\underline{i}s',\underline{i}t}=\xi_{\underline{i}s,\underline{i}t}.
\]
By induction, we get the desired result. \hfill$\square$
\end{pf}

\begin{remark} Note that the coefficients $a_{\underline{i},\underline{j}}$ and
$a_{\underline{i},\underline{k}}$ in the proof of
Lemma~\ref{L:generator} are smaller than $r!$. So
Lemma~\ref{L:generator} holds over $K$ whose characteristic equals
to $0$ or is greater than $r!$. \hfill$\square$
\end{remark}

\begin{Thm}
{\rm (i)} The natural action of $\mathcal{U}(\mathfrak{g}_{n})$ on
$\widetilde{E}^{\otimes r}$ induces a surjection
$\widetilde{\pi}:\mathcal{U}(\mathfrak{g}_{n})\rightarrow
\widetilde{S}(n,r)$.

{\rm (ii)} On $\mathcal{U}(\mathfrak{g}'_{n})$, we have
$\widetilde{\pi}^{r}=\widetilde{\rm{det}}_{a}^{\#}\circ
\widetilde{\pi}^{n+r}$.  Moreover, the following diagram commutes:
\[
{\setlength{\unitlength}{0.7pt}
\begin{picture}(480,200)
\put(50,130){$S(n,n+r)$}

\put(280,130){ $S(n,r) $}

\put(50,0){$\widetilde{S}(n,n+r)$}

\put(280,0){ $\widetilde{S}(n,r)$}

\put (170,30){$\mathcal{U}(\mathfrak{g}'_{n})$}

\put(170,160){$\mathcal{U}(\mathfrak{sl}_{n})$}

\put(125,134){ $\line(1,0){145}$}

\drawline(277,134)(272,139) \drawline(277,134)(272,129)
\drawline(272,134)(267,139) \drawline(272,134)(267,129)

\put(125,4){ $\line(1,0){145}$}

\drawline(277,4)(272,9) \drawline(277,4)(272,-1)
\drawline(272,4)(267,9) \drawline(272,4)(267,-1)

\put(85,115){$\line(0,-1){95}$}

\drawline(85,115)(80,110) \drawline(85,115)(90,110)
\drawline(85,110)(80,105) \drawline(85,110)(90,105)

\put(300,115){$\line(0,-1){95}$}

\drawline(300,115)(295,110) \drawline(300,115)(305,110)
\drawline(300,110)(295,105) \drawline(300,110)(305,105)

\drawline(100,146)(160,159)

 \drawline(100,146)(104,151) \drawline(100,146)(105,143)
 \drawline(105,147)(109,152) \drawline(105,147)(110,144)

\drawline(215,157)(280,146)

\drawline(280,146)(277,152) \drawline(280,146)(275,143)
\drawline(275,147)(272,153) \drawline(275,147)(270,144)

\dashline{5}(100,16)(165,29)

 \drawline(100,16)(104,21) \drawline(100,16)(105,13)

\dashline{5}(210,30)(280,16)

\drawline(280,16)(277,22) \drawline(280,16)(275,13)

\dashline{5}(185,150)(185,45)

\drawline(185,150)(180,145) \drawline(185,150)(190,145)
\drawline(185,145)(180,140) \drawline(185,145)(190,140)

\put(370,160){$\mathcal{U}({\mathfrak{gl}}_{n})$}

\drawline(365,159)(310,146)
\drawline(310,146)(314,151)\drawline(310,146)(315,143)
\drawline(315,147)(319,152) \drawline(315,147)(320,144)

\drawline(230,165)(355,165) \put(235,167.5){$\oval(20,5)[bl]$}
\put(235,171){$\oval(20,5)[tl]$} \drawline(355,165)(350,170)
\drawline(355,165)(350,160)

\put(370,30){$\mathcal{U}(\mathfrak{g}_{n})$}

\drawline(365,29)(310,16)
\drawline(310,16)(314,21)\drawline(310,16)(315,13)
\drawline(315,17)(319,22) \drawline(315,17)(320,14)

\dashline{5}(235,35)(355,35) \put(240,37.5){$\oval(20,5)[bl]$}
\put(240,41){$\oval(20,5)[tl]$} \drawline(355,35)(350,40)
\drawline(355,35)(350,30)

\put(385,145){$\line(0,-1){95}$}
 \drawline(385.5,145)(380.5,140) \drawline(385.5,145)(389.5,140)
 \drawline(385.5,140)(380.5,135) \drawline(385.5,140)(389.5,135)

\put(90,60){$\psi_{a}$} \put(195,90){$\eta_{a}$}
\put(305,60){$\psi_{a}$} \put(395,90){$\eta_{a}$}

\put(150,120){$\rm{det}^{*}$}
\put(150,-15){$\widetilde{\rm{det}}_{a}^{\#}$}
\put(130,28){$\widetilde{\pi}^{n+r}$}
\put(222,12){$\widetilde{\pi}^{r}$} \put(130,160){$\pi^{n+r}$}
\put(222,142){$\pi^{r}$}
\end{picture}}\]
\vskip3pt \noindent where $\eta_{a}$ is induced from the map denoted
by the same symbol defined in Section~\ref{S:semigroups}.

{\rm (iii)} When $r<n$, the restriction map
$\widetilde{\pi}=\widetilde{\pi}|_{\mathcal{U}(\mathfrak{g}'_{n})}:\mathcal{U}(\mathfrak{g}'_{n})\rightarrow
\widetilde{S}(n,r)$ is surjective.
\end{Thm}

\begin{pf} (i) It follows by Lemma~\ref{L:generator} (i).

(ii) To prove the first assertion it suffices to prove the two maps
coincide on a set of generators of $\mathcal{U}(\mathfrak{g}'_{n})$.
For $s=1,\ldots,n$, $t=s\pm 1$,
\[ \begin{array}{rcl}
\widetilde{\text{det}}_{a}^{\#}\circ\widetilde{\pi}^{n+r}(E_{st})&=&\widetilde{\text{det}}_{a}^{\#}(\sum_{\underline{p}\in
I(n,n+r-1)/\Sigma_{n+r-1}}\xi_{\underline{p}s,\underline{p}t})\\
&=&\sum_{\underline{p}\in
I(n,n+r-1)/\Sigma_{n+r-1}}\widetilde{\text{det}}_{a}^{\#}(\xi_{\underline{p}s,\underline{p}t}).
\end{array}
\]

If $\{p_{1},\ldots,p_{n+r-1}\}\nsupseteq\{1,\ldots,n\}$, then
$\widetilde{\text{det}}_{a}^{\#}(\xi_{\underline{p}s,\underline{p}t})=0$.
If $\underline{p}\sim_{\Sigma_{n+r-1}}(1\ldots n)\underline{q}$ for
some $\underline{q}\in I(n,r-1)$, then
$\widetilde{\text{det}}_{a}^{\#}(\xi_{\underline{p}s,\underline{p}t})=\xi_{\underline{q}s,\underline{q}t}$.
Hence
\[
\widetilde{\text{det}}_{a}^{\#}\circ\widetilde{\pi}^{n+r}(E_{st})=\sum_{\underline{q}\in
I(n,r-1)/\Sigma_{r-1}}\xi_{\underline{q}s,\underline{q}t}
=\widetilde{\pi}^{r}(E_{st}).
\]

Concerning the commutativity of the diagram, by Section~\ref{S:gln}
and Theorem~\ref{T:gpandalg} it remains to prove $\psi_{a}\circ
\widetilde{\pi}^{n+r}(E_{st})=\pi^{n+r}\circ\eta_{a}(E_{st})$ for
$s=1,\ldots,n$, $t\in\mathbb{Z}$. A direct check completes the
proof: if $s\neq t$, then by Lemma~\ref{L:imageofgen},
Proposition~\ref{P:psi} (ii), Lemma~\ref{L:imageofgen-gln} and the
definition of $\eta_{a}$ we have
\[\begin{array}{rcl}
\psi_{a}\circ\widetilde{\pi}^{n+r}(E_{st})
\hspace{-7pt}&=\hspace{-7pt}&\psi_{a}(\sum_{\underline{p}\in
I(n,n+r-1)/\Sigma_{n+r-1}}\xi_{\underline{p}s,\underline{p}t})\\
\hspace{-7pt}&=\hspace{-7pt}&\sum_{\underline{p}\in
I(n,n+r-1)/\Sigma_{n+r-1}}a^{\frac{t-\overline{t}}{n}}\xi_{\underline{p}s,\underline{p}\overline{t}}\\
\hspace{-7pt}&=\hspace{-7pt}&\pi^{n+r}(a^{\frac{t-\overline{t}}{n}}E_{s,\overline{t}})=\pi^{n+r}\circ\eta_{a}(E_{st});
\end{array}\]
if $s=t$, then by Lemma~\ref{L:imageofgen}
$\widetilde{\pi}^{n+r}(E_{ss})$ lies in $S(n,r)$, and hence by
Proposition~\ref{P:psi} (iii) and Lemma~\ref{L:eta} (iii) we have
$\psi_{a}\circ\widetilde{\pi}^{n+r}(E_{ss})=\widetilde{\pi}^{n+r}(E_{ss})=\pi^{n+r}(E_{ss})=\pi^{n+r}\circ\eta_{a}(E_{ss})$.

(iii) It follows from Lemma~\ref{L:generator} (ii). \hfill$\square$
\end{pf}

As a consequence, we have

\begin{Thm}
There is a surjective algebra homomorphism from
$\mathcal{U}(\widehat{\mathfrak{gl}}_{n})$ to
$\widetilde{S}(n,r)$. When $r<n$, the restriction of this
homomorphism to $\mathcal{U}(\widehat{\mathfrak{sl}}_{n})$ is also
surjective.
\end{Thm}

\begin{remark}
G.Lusztig~\cite{L} showed that when $r\geq n$ the quantized
$\widetilde{\pi}$ is not surjective.\hfill$\square$
\end{remark}

\section{Appendix 1 : Formal coalgebras}\label{S:formalcoalgebra}

Let $K$ be a field. In this section, when we say $K$-vector space we
mean a $K$-vector space $V$ with a fixed basis $\{v_{i}|i\in I\}$.
Let $\overline V$ denote the closure of $V$ with respect to formal
sums, i.e. $\overline{V}=\{\sum_{i\in I}{\lambda_{i}v_{i}}|
\lambda_{i}\in K\}$, and $i_{V}$ the canonical embedding from $V$ to
$\overline{V}$. For example, the closure of the polynomial ring
$K[X]$ with basis $\{X^{i}\}_{i\geq 0}$ is $K[[X]]$, the ring of
formal power series. Another example is, if $V$ is a finite
dimensional space, then $\overline{V}=V$.

 Denote by $V^{*}$ the dual of $V$, i.e.
$V^{*}=\{K$-linear maps $f:V\rightarrow K\}$, and by $V^{\#}$,
called {\em the dual of $V$ with finite support}, the subspace of
$V^{*}$ with basis $\{v_{i}^{*}|i\in I\}$ where
$v_{i}^{*}(v_{j})=\delta_{ij}$. Then $V^{*}=\overline{V^{\#}}$. Note
that $V^{*}=V^{\#}$ if $V$ is finite dimensional.

Let $V$, $W$ be two $K$-vector spaces with bases $\{v_{i}|i\in I\}$
and $\{w_{j}|j\in J\}$ respectively. A $K$-linear map
$f:V\rightarrow \overline{W}$ is called a {\em formal map} from $V$
to $W$. Such a map corresponds to a $J\times I$ matrix with entries
in $K$, say $M$. We say $f$ {\em is row finite} if there are only
finitely many nonzero entries in each row of $M$; we say $f$ {\em is
column finite} if there are only finitely many nonzero entries in
each column of $M$. In fact,  a formal map $f$ is column finite
 means exactly the image $\text{Im}(f)$ of $f$ lies in $W$. A linear
function $f\in V^{*}$ is row finite if and only if $f\in V^{\#}$.

The tensor product of $V$ and $W$ over $K$, denoted by $V\otimes W$,
has a natural basis $\{v_{i}\otimes w_{j}|i\in I, j\in J\}$. The
space $\overline{V} \otimes W$ is considered as a subspace of
$\overline{V} \otimes \overline{W}$ via $(\sum_{i\in
I}{\mu_{i}v_{i}})\otimes w_{j}\mapsto \sum_{i\in
I}{\mu_{i}v_{i}\otimes w_{j}}$, and the space $\overline{V} \otimes
\overline{W}$ is considered as a subspace of $\overline{V \otimes
W}$ via $(\sum_{i\in I}{\mu_{i}v_{i}})\otimes (\sum_{j\in
J}{\lambda_{j}w_{j}})\mapsto \sum_{i\in I,j\in
J}{\mu_{i}\lambda_{j}v_{i}\otimes w_{j}}$.

\begin{Prop}\label{P:tensor} We have a canonical isomorphism
$V^{\#}\otimes W^{\#}\rightarrow (V\otimes W)^{\#}$,
$v_{i}^{*}\otimes w_{j}^{*}\mapsto (v_{i}\otimes w_{j})^{*}$. We
identify these two spaces.  \hfill$\square$
\end{Prop}

We have $V\cong (V^{\#})^{\#}$ canonically. Identifying these two
spaces, we have
$\overline{V}=(V^{\#})^{*}$. 
Thus, as a consequence of Proposition~\ref{P:tensor}, we have
$\overline{V\otimes V}=((V\otimes V)^{\#})^{*}= (V^{\#}\otimes
V^{\#})^{*}$.

Let $f:V\rightarrow \overline{W}$, $v_{i}\mapsto \sum_{j\in
J}{m_{ji}w_{j}}$, be a row finite formal map. Then we can extend $f$
to a linear map $\overline{f}:\overline{V}\rightarrow \overline{W}$,
$\sum_{i\in I}{\lambda_{i}v_{i}} \mapsto \sum_{j\in j}{(\sum_{i\in
I}{\lambda_{i}m_{ji}})w_{j}}$.

\begin{Thm}\label{T:formalmaps} Let $V$, $W$ be two $K$-vector spaces with basis
$\{v_{i}|i\in I\}$ and $\{w_{j}|j\in J\}$ respectively, and
$f:V\rightarrow \overline{W}$ a formal map. Then $f^{\#}:
W^{\#}\rightarrow V^{*}$, $\alpha\mapsto \overline{\alpha}\circ f$
is a formal map from $W^{\#}$ to $V^{\#}$. Moreover,

{\rm(i)} If \ $f$ is column finite, then $f^{\#}$ is row finite.

{\rm(ii)} If \ $f$ is row finite, then $f^{\#}$ is column finite.

{\rm(iii)} If \ ${\rm Im}(f)\supseteq W$, then $f^{\#}$ is
injective.

{\rm (iv)} If \ $f$ is row finite and $\overline{f}$ is injective,
then ${\rm Im}(f^{\#})=V^{\#}$.

Let $g:W\rightarrow\overline{U}$ be a row finite formal map.

{\rm(v)} If \ $f$ is also row finite, then $\overline{g}\circ f$ is
row finite and $\overline{\overline{g}\circ
f}=\overline{g}\circ\overline{f}:\overline{V}\rightarrow
\overline{U}$.

{\rm(vi)} We have $(\overline{g}\circ f)^{\#}=f^{\#}\circ
g^{\#}:U^{\#}\rightarrow V^{*}$.

\end{Thm}

\begin{pf}

Let $M=(m_{ji})_{i\in I,j\in J}$ be the matrix corresponding to
$f$. Then
\[
f^{\#}(w_{j}^{*})(v_{i})=w_{j}^{*}(f(v_{i}))=w_{j}^{*}(\sum_{l\in
J}{m_{li}w_{l}})=m_{ji}.
\]
So $f^{\#}(w_{j}^{*})=\sum_{i\in I}{m_{ji}v_{i}^{*}}$, and hence
$f^{\#}$ corresponds to $M^{tr}$, the transpose of $M$. (i) and (ii)
follow immediately.

(iii) Assume that $\text{Im}(f)\supseteq W$. Let $\alpha\in W^{\#}$
such that $f^{\#}(\alpha)=0$, i.e. $f^{\#}(\alpha)(v)=0$, for any
$v\in V$. Therefore $\alpha(w_{j})=0$ for any $j\in J$. So
$\alpha=0$, i.e. $f^{\#}$ is injective.

(iv) Assume $f$ is row finite and $\overline{f}$ is injective. It
follows by (ii) that $\text{Im}(f^{\#})\subseteq V^{\#}$. Suppose
$\text{Im}(f)$ is a proper subspace of $V^{\#}$. Then there exists a
nonzero element $v\in \overline{V}=(V^{\#})^{*}$ such that
$\overline{f^{\#}(\alpha)}(v)=\overline{\overline{\alpha}\circ
f}(v)=0$ for any $\alpha\in W^{\#}$. Therefore by (v) we have
$\overline{\alpha}(\overline{f}(v))=0$ for any $\alpha\in W^{\#}$,
and hence $\overline{f}(v)=0$, contradicting the injectivity of
$\overline{f}$.

(v) We have
\[\begin{array}{c}\overline{g}\circ f(v_{i})=\overline{g}(\sum_{j\in
J}m_{ji}w_{j})=\sum_{j\in J}m_{ji}g(w_{j}).\end{array}\] Therefore
the matrix of $\overline{g}\circ f$ is the product two row finite
matrices: the matrices of $g$ and $f$, and is again row finite. The
proof for the equality is straightforward:
\[\begin{array}{rcl}
\overline{\overline{g}\circ f}(\sum_{i\in
I}\lambda_{i}v_{i})\hspace{-7pt}&=\hspace{-7pt}&\sum_{i\in
I}\lambda_{i}\overline{g}\circ
f(v_{i})=\sum_{i\in I,j\in J}\lambda_{i}m_{ji}g(w_{j})\\
\hspace{-7pt}&=\hspace{-7pt}&\overline{g}(\sum_{i\in I,j\in
J}\lambda_{i}m_{ji}w_{j})
=\overline{g}(\sum_{i\in I}\lambda_{i}f(v_{i}))\\
\hspace{-7pt}&=\hspace{-7pt}&\overline{g}\circ\overline{f}(\sum_{i\in
I}\lambda_{i}v_{i}).
\end{array}
\]

(vi) It follows by (ii) that $\text{Im}(g^{\#})\subseteq W^{\#}$,
and hence $f^{\#}\circ g^{\#}$ is well-defined. For $\beta\in
U^{\#}$, we have
\[ \begin{array}{c}
 f^{\#}\circ g^{\#}(\beta)
=f^{\#}(\overline{\beta}\circ g)=\overline{\overline{\beta}\circ
g}\circ f =\overline{\beta}\circ\overline{g}\circ f
=(\overline{g}\circ
 f)^{\#}(\beta).\end{array}
\]
Hence $(\overline{g}\circ f)^{\#}=\overline{f^{\#}}\circ g^{\#}$.
\hfill$\square$
\end{pf}

A $K$-linear map $f:V\rightarrow W$ can be considered as a formal
map from $V$ to $W$ via the embedding $i_{W}:W\rightarrow
\overline{W}$.

\begin{Cor}\label{C:formalmaps}
Let $f: V\rightarrow W$ be a row finite $K$-linear map. Then
$f^{\#}: W^{\#}\rightarrow V^{\#}$, $\alpha\mapsto \alpha\circ f$,
is a row finite $K$-linear map. Moreover, if \ $\overline{f}$ is
injective, then $f^{\#}$ is surjective; if \ $f$ is surjective, then
$f^{\#}$ is injective; if $g:W\rightarrow U$ is a row finite
$K$-linear map, then $(g\circ f)^{\#}=f^{\#}\circ g^{\#}$.
\end{Cor}
\begin{pf} The map $f: V\rightarrow W$ is a row finite $K$-linear
map, and hence is a row finite and column finite formal map from $V$
to $W$. Applying Theorem~\ref{T:formalmaps} we get the desired
results. \hfill$\square$
\end{pf}

\begin{Def}\label{D:formalcoalgebra}
A {\em formal $K$-coalgebra} is a $K$-linear space $A$ together with
two formal maps $\Delta:A\rightarrow \overline{A\otimes A}$ and
$\epsilon:A\rightarrow K$ such that

{\rm(i)} the maps $\Delta$ and $\epsilon$ are row finite,

{\rm(ii)} the following diagrams commute:
\[{\setlength{\unitlength}{0.7pt}
\begin{picture}(180,90)
\put(17,60){$A$} \put(100,60){ $\overline{A\otimes A} $}
\put(0,0){$\overline{A\otimes A} $} \put(100,0){ $\overline{A\otimes
A\otimes A}$} \put(30,30){\footnotesize$\Delta$}
\put(60,67){\footnotesize$\Delta$}
\put(128,30){\footnotesize$\overline{\Delta\otimes \rm{id}}$}
\put(52,7){\footnotesize$\overline{\rm{id}\otimes \Delta}$}
\drawline(50,5)(95,5) \drawline(95,5)(91,1) \drawline(95,5)(91,9)
\drawline(40,65)(95,65) \drawline(95,65)(91,61)
\drawline(95,65)(91,69) \drawline(23,50)(23,20)
\drawline(23,20)(19,24) \drawline(23,20)(27,24)
\drawline(125,50)(125,20) \drawline(125,20)(121,24)
\drawline(125,20)(129,24)
\end{picture} \hskip5mm\begin{picture}(130,90)
\put(17,60){$A$} \put(17,0){$\overline{A}$} \put(75,60){
$\overline{A\otimes A} $} \put(75,0){ $\overline{A\otimes K}$}
\put(48,7){\footnotesize$\cong$} \put(48,67){\footnotesize$\Delta$}
\put(28,30){\footnotesize$i_{A}$} \drawline(35,65)(70,65)
\drawline(70,65)(66,61) \drawline(70,65)(66,69)
\put(25,50){$\oval(5,10)[t]$} \drawline(22.5,50)(22.5,20)
\drawline(22.5,20)(18.5,24) \drawline(22.5,20)(26.5,24)
\drawline(35,5)(70,5) \drawline(70,5)(66,1) \drawline(70,5)(66,9)
\drawline(100,50)(100,20) \drawline(100,20)(96,24)
\drawline(100,20)(104,24)
\put(105,30){\footnotesize$\overline{\rm{id}\otimes \epsilon}$}
\end{picture}\hskip7mm\begin{picture}(130,90)
\put(17,60){$A$} \put(17,0){$\overline{A}$} \put(75,60){
$\overline{A\otimes A} $} \put(75,0){ $\overline{K\otimes A}$}
\put(48,7){\footnotesize$\cong$} \put(48,67){\footnotesize$\Delta$}
\drawline(35,65)(70,65) \drawline(70,65)(66,61)
\drawline(70,65)(66,69) \put(25,50){$\oval(5,10)[t]$}
\drawline(22.5,50)(22.5,20) \drawline(22.5,20)(18.5,24)
\drawline(22.5,20)(26.5,24) \drawline(35,5)(70,5)
\drawline(70,5)(66,1) \drawline(70,5)(66,9)
\drawline(100,50)(100,20) \drawline(100,20)(96,24)
\drawline(100,20)(104,24) \put(28,30){\footnotesize$i_{A}$}
\put(105,30){\footnotesize$\overline{\epsilon\otimes \rm{id}}$}
\end{picture}}\]
We call $\Delta$ comultiplication, and $\epsilon$ counit.

\end{Def}

\begin{remark} (i) A $K$-coalgebra with row finite comultiplication
and row finite counit is a formal $K$-coalgebra. A formal
$K$-coalgebra with column finite comultiplication is a
$K$-coalgebra.

(ii) Let $A$ be a formal $K$-coalgebra with comultiplication
$\Delta$ and counit $\epsilon$. Then $\Delta: A\rightarrow
\overline{A\otimes A}$ can be lifted to $\overline{\Delta}:
\overline{A}\rightarrow \overline{A\otimes A}$. The map
$\text{id}\otimes\Delta: A\otimes A\rightarrow A\otimes
\overline{A\otimes A}$ can be lifted to $\overline{A\otimes A\otimes
A}$. The same holds for $\Delta\otimes \text{id}:A\otimes
A\rightarrow \overline{A\otimes A}\otimes A$. \hfill$\square$
\end{remark}

\begin{Def}\label{D:formalhomomorphism}
Let $A$, $B$ be two formal $K$-coalgebras with comultiplications
$\Delta_{A}$, $\Delta_{B}$ and counits $\epsilon_{A}$,
$\epsilon_{B}$ respectively. A {\em formal homomorphism of formal
$K$-coalgebras} from $A$ to $B$ is a formal map $f:A\rightarrow
\overline{B}$ such that

{\rm(i)} the map $f$ is row finite,

{\rm(ii)} the following diagrams commute:
\[
{\setlength{\unitlength}{0.7pt}
\begin{picture}(160,90)
\put(17,60){$A$}

\put(115,60){ $\overline{B}$}

\put(0,0){$\overline{A\otimes A} $}

\put(100,0){ $\overline{B\otimes B}$}

\put(28,30){\footnotesize$\Delta_{A}$}

\put(70,67){\footnotesize$f$}

\put(128,30){\footnotesize$\overline{\Delta}_{B}$}

\put(52,7){\footnotesize$\overline{f\otimes f}$}

\drawline(40,65)(105,65) \drawline(105,65)(101,61)
\drawline(105,65)(101,69)

\drawline(50,5)(95,5) \drawline(95,5)(91,1) \drawline(95,5)(91,9)

\drawline(23,50)(23,20) \drawline(23,20)(19,24)
\drawline(23,20)(27,24)

\drawline(125,50)(125,20) \drawline(125,20)(121,24)
\drawline(125,20)(129,24)

\end{picture}\ \ \ \  \begin{picture}(130,90)
\put(17,60){$A$}

\put(17,0){$K$}

\put(75,60){ $\overline{B}$}

\put(75,0){ $K$}

\put(48,7){\footnotesize$\rm{id}$}

\put(48,67){\footnotesize$f$}

\drawline(35,65)(70,65) \drawline(70,65)(66,61)
\drawline(70,65)(66,69)

\drawline(22.5,50)(22.5,20) \drawline(22.5,20)(18.5,24)
\drawline(22.5,20)(26.5,24)

\drawline(35,5)(70,5) \drawline(70,5)(66,1) \drawline(70,5)(66,9)

\drawline(85,50)(85,20) \drawline(85,20)(81,24)
\drawline(85,20)(89,24)

\put(28,30){$\epsilon_{A}$} \put(92,30){$\overline{\epsilon}_{B}$}
\end{picture}}
\]
If \ $f$ is column finite, then $f: A\rightarrow B$ is called a {\em
homomorphism of formal $K$-coalgebras} from $A$ to $B$.

\end{Def}

\begin{Thm}\label{T:fromcoalgebratoalgebra} Let $A$, $B$ be two formal
$K$-coalgebras with comultiplications $\Delta_{A}$, $\Delta_{B}$
and counits $\epsilon_{A}$, $\epsilon_{B}$ respectively.

{\rm(i)} The dual $A^{\#}$ of $A$ with finite support is a
$K$-algebra with multiplication $\Delta_{A}^{\#}$ and unit
$\epsilon_{A}^{\#}$. Precisely, $\Delta_{A}^{\#}(f\otimes
f')(a)=(f\otimes f')(\Delta_{A}(a))=\sum f(a_{1})f'(a_{2})$, where
$f,f'\in A^{\#}$, $a\in A$, and $\Delta_{A}(a)=\sum a_{1}\otimes
a_{2}$ is the Sweedler symbol.

{\rm(ii)} Let $f:A\rightarrow \overline{B}$ be a formal homomorphism
of formal $K$-coalgebras from $A$ to $B$. Then $f^{\#}:
B^{\#}\rightarrow A^{\#}$ is a homomorphism of $K$-algebras.
Moreover, if \ $\overline{f}$ is injective, then $f^{\#}$ is
surjective; if \ ${\rm Im}(f)\supseteq B$, then $f^{\#}$ is
injective.

{\rm(iii)} Let $f:A\rightarrow B$ be a homomorphism of formal
$K$-coalgebras. Then $f^{\#}: B^{\#}\rightarrow A^{\#}$ is a row
finite homomorphism of $K$-algebras. Moreover, if \ $\overline{f}$
is injective, then $f^{\#}$ is surjective; if \ $f$ is surjective,
then $f^{\#}$ is injective.
\end{Thm}
\begin{pf}
(i) Since $\Delta_{A}:A\rightarrow\overline{A\otimes A}$ and
$\epsilon_{A}:A\rightarrow K$ are row finite, it follows by
Theorem~\ref{T:formalmaps} (ii) that $\Delta_{A}^{\#}:(A\otimes
A)^{\#}=A^{\#}\otimes A^{\#}\rightarrow \overline{A^{\#}}$ and
$\epsilon_{A}^{\#}:K\rightarrow \overline{A^{\#}}$ are column
finite. Namely, $\Delta_{A}^{\#}:(A\otimes A)^{\#}=A^{\#}\otimes
A^{\#}\rightarrow A^{\#}$ and $\epsilon_{A}^{\#}:K\rightarrow
A^{\#}$ are two $K$-linear maps. Now the dual of the diagrams in
Definition~\ref{D:formalcoalgebra} (ii) consists of the axioms for
$A^{\#}$ to be an algebra.

(ii) Since $f$ is row finite, it follows by
Theorem~\ref{T:formalmaps} (ii) that $f^{\#}:B^{\#}\rightarrow
\overline{A^{\#}}$ is column finite, i.e. $f^{\#}:B^{\#}\rightarrow
A^{\#}$ is a $K$-linear map. Further, the dual of the diagrams in
Definition~\ref{D:formalhomomorphism} (ii) consists of the axioms
for $f^{\#}$ to be an algebra homomorphism. The last two assertions
follow from Theorem~\ref{T:formalmaps} (iii) (iv).

(iii) It follows from (ii) that $f^{\#}$ is an algebra homomorphism,
and it follows from Theorem~\ref{T:formalmaps} (i) that $f^{\#}$ is
row finite. The last two assertions follow from
Corollary~\ref{C:formalmaps}. \hfill $\square$

\end{pf}

\section{Appendix 2 : Mackey formula for certain endomorphism
algebras}\label{S:mackey}

This appendix contains a generalization of \cite{JAG2} (4.11).

Let $G$ be a group. Let $I$ be a set with a right $G$-action such
that the stabilizer of any $i\in I$ is a finite subgroup of $G$.

Let $K$ be a field and $V$ a $K$-vector space with basis
$\{v_{i}\}_{i\in I}$. Then $V$ is right $G$-module over $K$ via
$v_{i}^{g}=v_{ig}$ where $g\in G$, and $i\in I$.  Let
$B=\text{End}_{K}(V)$. For $i$,$j\in I$, let $x_{ij}\in B$ be the
$K$-map defined by $x_{ij}(v_{l})=\delta_{jl}v_{i}$ where $l\in I$.
Then
\[
\begin{array}{c}
B=\{\sum_{i,j\in I}\lambda_{ij}x_{ij}|{\rm\ for\ fixed\ } j,
\lambda_{ij}=0 {\rm\ for\ almost\ all\ } i \}
\end{array}
\]

The right $G$-module structure on $V$ induces a right $G$-module
structure on $B$ satisfying $(bb')^{g}=b^{g}(b')^{g}$ for any
$b,b'\in B$ and $g\in G$. Precisely, $(\sum_{i,j\in
I}\lambda_{ij}x_{ij})^{g}=\sum_{i,j\in I}\lambda_{ij}x_{ig,jg}$. Let
$A$ be the subring of $B$ with basis $\{x_{ij}\}_{i,j\in I}$. It is
easy to see that $A$ is a $G$-submodule of $B$.

For a subgroup $H\leq G$, let
\[
\begin{array}{c}
A_{H}=\{a\in A|a^{h}=a, \forall h\in H\}\\
B_{H}=\{b\in B|b^{h}=b, \forall h\in H\}.
\end{array}
\]

For a subgroup chain $H_{1}\leq H_{2}\leq G$ of $G$, let
\[
T_{H_{1},H_{2}} : A_{H_{1}}\rightarrow B_{H_{2}},~~ a\mapsto
\sum_{g\in H_{1}\backslash H_{2}}a^{g}.
\]
The map $T_{H_{1},H_{2}}$ is well-defined because of the following.
Let $a=\sum_{i,j}\lambda_{ij}x_{ij}\in A_{H_{1}}$, then
\[
\begin{array}{rcl}
T_{H_{1},H_{2}}(a)
\hspace{-7pt}&=\hspace{-7pt}&\sum_{g}(\sum_{i,j}\lambda_{ij}x_{ij})^{g}=\sum_{i,j}\lambda_{ij}\sum_{g}(x_{ij})^{g}\\
\hspace{-7pt}&=\hspace{-7pt}&\sum_{i,j,g}\lambda_{ij}x_{ig,jg}
=\sum_{p,q}(\sum_{i,j,g:ig=p,jg=q}\lambda_{ij})x_{pq}.
\end{array}
\]
But for fixed $(i,j)$ and $(p,q)$ there are only finitely many $g$
such that $(ig,jg)=(p,q)$.

\begin{Lem}\label{L:move} Let $H_{1}\leq H_{2}\leq G$ be a subgroup chain of $G$
and $a\in A_{H_{1}}$, $b\in B_{H_{2}}$.

{\rm (i)} If \ $ab\in A_{H_{1}}$, then
$T_{H_{1},H_{2}}(ab)=T_{H_{1},H_{2}}(a)b$.

{\rm (ii)} If \ $ba\in A_{H_{1}}$, then
$T_{H_{1},H_{2}}(ba)=bT_{H_{1},H_{2}}(a)$.

\begin{pf} We only prove (i); the proof for (ii) is similar. Assume
$ab\in A_{H_{1}}$, then
\[
\begin{array}{rcl}
T_{H_{1},H_{2}}(ab) \hspace{-7pt}&=\hspace{-7pt}&\sum_{g\in
H_{1}\backslash H_{2}}(ab)^{g}=\sum_{g\in H_{1}\backslash
H_{2}}(a^{g}b^{g})\\
\hspace{-7pt}&=\hspace{-7pt}&\sum_{g\in H_{1}\backslash
H_{2}}(a^{g}b)=(\sum_{g\in
H_{1}\backslash H_{2}}a^{g})b\\
\hspace{-7pt}&=\hspace{-7pt}&T_{H_{1},H_{2}}(a)b.\end{array}
\]
\hfill$\square$
\end{pf}

\end{Lem}

\begin{Lem}\label{L:transitivity}
Let $H_{1}\leq H_{2}\leq H_{3}\leq G$ be a subgroup chain of $G$,
and $a\in A_{H_{1}}$. Assume $T_{H_{1},H_{2}}(a)\in A_{H_{2}}$,
then $T_{H_{2},H_{3}}(T_{H_{1},H_{2}}(a))=T_{H_{1},H_{3}}(a)$.
\end{Lem}
\begin{pf}
Let $C_{1}$ be a set of representatives of cosets $H_{1}\backslash
H_{2}$, and $C_{2}$ be a set of representatives of cosets
$H_{2}\backslash H_{3}$, then $C=\{c_{1}c_{2}|c_{1}\in
C_{1},c_{2}\in C_{2}\}$ is a set of representatives of cosets
$H_{1}\backslash H_{3}$. This finishes the proof. \hfill$\square$
\end{pf}

\begin{Lem}\label{L:compare}
Let $H_{1}\leq H_{3}\leq G$ and $H_{2}\leq H_{3}\leq G$ be two
subgroup chains of $G$, and $a\in A_{H_{1}}$. Then
\[
T_{H_{1},H_{3}}(a)=\sum_{w\in H_{1}\backslash
H_{3}/H_{2}}T_{H_{1}^{w}\cap H_{2},H_{2}}(a^{w}).
\]
\end{Lem}
\begin{pf} Let $W$ be a set of representatives of $H_{1}\backslash
H_{3}/H_{2}$. For each $w\in W$, let $C_{w}$ be a set of
representatives of $(H_{1}^{w}\cap H_{2})\backslash H_{2}$. Then
$\bigcup_{w\in W} wC_{w}$ is a set of representatives of
$H_{1}\backslash H_{3}$. Indeed, for $h\in H_{3}$ there exists
$h_{1}\in H_{1}$, $w\in W$ and $h_{2}\in H_{2}$ such that
$h=h_{1}wh_{2}$. Since $h_{2}\in H_{2}$, there exists $g\in
H_{1}^{w}\cap H_{2}$ and $c\in C_{w}$ such that $h_{2}=gc$. Thus
$h=h_{1}wgc=h_{1}wgw^{-1}wc\in H_{1}wc$. Conversely, suppose
$wc=w'c'$ where $w,w'\in W$, $c\in C_{w}$ and $c'\in C_{w'}$. Then
$w=w'c'c^{-1}$. Since $W$ is a set of representatives of
$H_{1}\backslash H_{3}/H_{2}$, and $c'c^{-1}\in H_{2}$, we have
$w=w'$ and $c=c'$.\hfill$\square$
\end{pf}

\begin{Lem}\label{L:mackey} (Mackey's formula)
Let $H_{1}\leq H_{3}\leq G$ and $H_{2}\leq H_{3}\leq G$ be two
subgroups chains of $G$, and $a\in A_{H_{1}}$, $b\in A_{H_{2}}$.
If \ $aT_{H_{2},H_{3}}(b)\in A_{H_{1}}$, $T_{H_{1}\cap
H_{2}^{w},H_{1}}(ab^{w})=0$ for almost all $w$, and $T_{H_{1}\cap
H_{2}^{w},H_{1}}(ab^{w})\in A_{H_{1}}$ for all $w$, then
\[
T_{H_{1},H_{3}}(a)T_{H_{2},H_{3}}(b)=\sum_{w\in H_{1}\backslash
H_{3}/H_{2}}T_{H_{1}\cap H_{2}^{w},H_{3}}(ab^{w}).
\]

\end{Lem}

\begin{pf} First note that $b\in
A_{H_{2}}$ implies $b^{w}\in A_{H_{2}^{w}}$. Thus $a$ and $b^{w}$
are in $A_{H_{1}\cap H_{2}^{w}}$ and so is their product $ab^{w}$.

By Lemma~\ref{L:compare} and Lemma~\ref{L:move} we have
\[\begin{array}{rcl}
aT_{H_{2},H_{3}}(b)&=&a\sum_{w\in H_{2}\backslash
H_{3}/H_{1}}T_{H_{1}\cap H_{2}^{w},H_{1}}(b^{w})\\
&=&\sum_{w\in H_{2}\backslash H_{3}/H_{1}}T_{H_{1}\cap
H_{2}^{w},H_{1}}(ab^{w}).
\end{array}\]

Since $aT_{H_{2},H_{3}}(b)\in A_{H_{1}}$, it follows that
\[\begin{array}{rcl}
T_{H_{1},H_{3}}(a)T_{H_{2},H_{3}}(b)
&=&T_{H_{1},H_{3}}(aT_{H_{2},H_{3}}(b))\\
&=&T_{H_{1},H_{3}}(\sum_{w\in H_{2}\backslash
H_{3}/H_{1}}T_{H_{1}\cap H_{2}^{w},H_{1}}(ab^{w}))\\
 &=&\sum_{w\in
H_{2}\backslash
H_{3}/H_{1}}T_{H_{1},H_{3}}(T_{H_{1}\cap H_{2}^{w},H_{1}}(ab^{w}))\\
&=&\sum_{w\in H_{2}\backslash H_{3}/H_{1}}T_{H_{1}\cap
H_{2}^{w},H_{3}}(ab^{w}),
\end{array}\]
where the last equality follows by Lemma~\ref{L:transitivity}.
\hfill$\square$
\end{pf}

\begin{ac}
The author would like to thank Steffen Koenig for his inspiring
and encouraging supervision.
\end{ac}

\bibliographystyle{amsplain}

\begin{thebibliography}{99}

\bibitem{CL} R.W.Carter and G.Lusztig, On the modular
representations of the general linear and symmetric groups, Math.
Zeit. 136 (1974), 193--242.

\bibitem{CP} V.Chari and A.Pressley, Quantum affine algebras
and affine Hecke algebras, Pacific J. Math. 174 (1996), no. 2,
295--326.


\bibitem{Du} J.Du, $q$-Schur algebras, asymptotic forms and quantum $SL_n$,
J. Algebra 177 (1995), 385--408.

\bibitem{DR} S.R.Doty and R.M.Green, Presenting affine $q$-Schur
algebras, arXiv:math.QA/0603002.

\bibitem{GRV} V.Ginzburg, N.Reshetikhin and E.Vasserot,
Quantum groups and flag varieties, Contemp. Math. 175 (1994),
101--130.

\bibitem{GV} V.Ginzburg and E.Vasserot, Langlands reciprocity
for affine quantum groups of type $A_{n}$, Internat. Math. Res.
Notices (1993) no.3, 67--85.

\bibitem{JAG} J.A.Green, Polynomial representations of $GL_n$, Lecture
notes on Mathematics 830, 1980.

\bibitem{JAG2} J.A.Green, Some remarks on defect groups,
Math. Zeit. (1968) 107, 133--150.

\bibitem{JAG3} J.A.Green, On certain subalgebras of the
Schur algebra, J. Algebra 131 (1990), 265--280.


\bibitem{JAG5} J.A.Green, Combinatorics and the Schur
algebra, Journal of Pure and Applied Algebra 88 (1993), 89--106.

\bibitem{RMG} R.M.Green, The affine $q$-Schur algebra, J. Algebra
215
(1999), no. 2, 379--411.

\bibitem{RMG3} R.M.Green, $q$-Schur algebras and quantized enveloping
algebras, Ph.D. thesis, University of Warwick, 1995.

\bibitem{RMG4} R.M.Green, Hyperoctahedral Schur algebras, J.Algebra
192 (1997), 418--438.

\bibitem{Ka} V.G.Kac, Infinite dimensional Lie algebras, 3rd edition,
Cambridge University Press, Cambridge, 1990.

\bibitem{SK} S.Koenig, Blocks of category O, double
centralizer properties, and Enright's completions, Proceedings of
NATO-ASI (Constanta, 2000), Algebras - Representation Theory,
113--134, Kluwer (2001).

\bibitem{Liu} Q.Liu, Schur algebras of classical groups II,
preprint.

\bibitem{L} G.Lusztig, Aperiodicity in quantum affine
$\mathfrak{gl}_{n}$, Sir Michael Atiyah: a great mathematician of
the twentieth century. Asian J. Math. 3 (1999) no.1, 147--177.

\bibitem{L2} G.Lusztig, Transfer maps for quantum affine
$\mathfrak{sl}_{n}$, Repesentations and quantizations (Shanghai,
1998), 341--356, China High. Educ. Press, Beijing (2000).

\bibitem{McGerty} K.McGerty, q-Schur algebras and quantum
Frobenius, Adv. Math., to appear.

\bibitem{VV} M.Varagnolo and E.Vasserot, On the decomposition matrices of
the quantized Schur algebra, Duke Math. J. 100 (1999), no.2,
267--297.


\end{thebibliography}

{Dong Yang}

{Department of Mathematical Sciences, Tsinghua University,
Beijing100084, P.R.China.}

{Department of Mathematics, University of Leicester, University
Road, Leicester LE1 7RH, United Kingdom.}

 {\it Email address : }
{yangdong98@mails.tsinghua.edu.cn.}

\end{document}